\documentclass[a4paper,11pt]{article}

\usepackage{amssymb}
\usepackage{amsthm}
\usepackage{bbm}
\usepackage{hyperref}
\usepackage{url}
\usepackage[all]{xy}

\title{Conjugacy $p$-separability of right-angled Artin groups and applications}
\author{Emmanuel Toinet}
\date{}

\makeatletter 
\newtheorem*{rep@theorem}{\rep@title} 
\newcommand{\newreptheorem}[2]{
\newenvironment{rep#1}[1]{
\def\rep@title{#2 \ref{##1}}
\begin{rep@theorem}}
{\end{rep@theorem}}} 
\makeatother  

\makeatletter 
\newcommand{\subjclass}[2][2010]{   
\let\@oldtitle\@title  
\gdef\@title{\@oldtitle\footnotetext{#1 
\emph{Mathematics Subject Classification.} #2.}}} 
\newcommand{\keywords}[1]{  
\let\@@oldtitle\@title 
\gdef\@title{\@@oldtitle\footnotetext{\emph{Key words and phrases.} #1.}}
} 
\makeatother 


\begin{document}

\maketitle

\begin{abstract}
We prove that every subnormal subgroup of $p$-power index in a right-angled Artin group is conjugacy $p$-separable. As an application, we prove that every right-angled Artin group is conjugacy separable in the class of torsion-free nilpotent groups. As another application, we prove that the outer automorphism group of a right-angled Artin group is virtually residually $p$-finite. We also prove that the Torelli group of a right-angled Artin group is residually torsion-free nilpotent, hence residually $p$-finite and bi-orderable.
\end{abstract}

\section{Introduction}

\newtheorem{definition}{Definition}[section]
\newtheorem{theorem}[definition]{Theorem}
\newreptheorem{theorem}{Theorem}
\newtheorem{proposition}[definition]{Proposition}
\newtheorem{lemma}[definition]{Lemma}
\newtheorem{corollary}[definition]{Corollary}
\newreptheorem{corollary}{Corollary}
\theoremstyle{remark}
\newtheorem{remark}[definition]{Remark}

\hspace{5mm}Let $\Gamma$ = ($V$,$E$) be a finite simplicial graph. The \textit{right-angled Artin group} associated to $\Gamma$ is the group $G_\Gamma$ defined by the presentation:
\begin{center}
$G_\Gamma$ = $\langle$ $V$ $\mid$ $vw$ = $wv$, $\forall$ \{$v$,$w$\} $\in$ $E$ $\rangle$.
\end{center}

Note that, if $\Gamma$ is a discrete graph, then $G_\Gamma$ is the free group $F_r$ on $r$ generators (where $r$ = $\left|V\right|$), and if $\Gamma$ is a complete graph, then $G_\Gamma$ is the free abelian group $\mathbbm{Z}^r$. Thus, right-angled Artin groups can be seen as interpolating between free groups and free abelian groups. The \textit{rank} of $G_\Gamma$ is by definition the number of vertices of $\Gamma$. A \textit{special subgroup} of $G_\Gamma$ is a subgroup generated by a subset $W$ of the set of vertices $V$ of $\Gamma$ -- it is naturally isomorphic to the right-angled Artin group $G_{\Gamma(W)}$, where $\Gamma(W)$ denotes the full subgraph of $\Gamma$ spanned by $W$. Let $v$ be a vertex of $\Gamma$. The \textit{link} of $v$, denoted by $link(v)$, is the subset of $V$ consisting of all vertices that are adjacent to $v$. The \textit{star} of $v$, denoted by $star(v)$, is $link(v)\cup\{v\}$. We refer to \cite{C} for a general survey of right-angled Artin groups.
\\

Little is known about the automorphim groups of right-angled Artin groups. In 1989, Servatius conjectured a generating set for $Aut(G_\Gamma)$ (see \cite{Ser}). He proved his conjecture in certain special cases, for example when the graph is a tree. Thereafter Laurence proved the conjecture in the general case (see \cite{L}). Charney and Vogtmann showed that $Out(G_\Gamma)$ is virtually torsion-free and has finite virtual cohomological dimension (see \cite{CV1}). Day gave a finite presentation for $Aut(G_\Gamma)$ (see \cite{D1}). More recently, Minasyan proved that $Out(G_\Gamma)$ is residually finite (see \cite{M}). This result was obtained independently by Charney and Vogtmann in \cite{CV2}, where they also proved that, for a large class of graphs, $Out(G_\Gamma)$ satisfies the Tits alternative.
\\

Let $\mathcal{K}$ be a class of group. A group $G$ is said to be \textit{residually $\mathcal{K}$} if for all $g$ $\in$ $G\setminus$\{1\}, there exists a homomorphism $\varphi$ from $G$ to some group of $\mathcal{K}$ such that $\varphi(g)$ $\neq$ 1. Note that if $\mathcal{K}$ is the class of all finite groups, this notion reduces to residual finiteness.

For a group $G$ and for $g$, $h$ $\in$ $G$, we use the notation $g$ $\sim$ $h$ to mean that $g$ and $h$ are conjugate. A group $G$ is said to be \textit{conjugacy $\mathcal{K}$-separable} (or \textit{conjugacy separable in the class $\mathcal{K}$}) if for all $g$, $h$ $\in$ $G$, either $g$ $\sim$ $h$, or there exists a homomorphism $\varphi$ from $G$ to some group of $\mathcal{K}$ such that $\varphi(g)$ $\nsim$ $\varphi(h)$. Note that if $\mathcal{K}$ is the class of all finite groups, this notion reduces to conjugacy separability. Clearly, if a group is conjugacy $\mathcal{K}$-separable, then it is residually $\mathcal{K}$.
\\

Our focus here is on conjugacy separability in the class of finite $p$-groups. Let $p$ be a prime number. If $\mathcal{K}$ is the class of all finite $p$-groups, then, instead of saying ``$G$ is residually $\mathcal{K}$", we shall say that $G$ is \textit{residually $p$-finite}. Note that this implies residually finite as well as residually nilpotent. Instead of saying ``$G$ is conjugacy $\mathcal{K}$-separable", we shall say that $G$ is \textit{conjugacy $p$-separable}. Following Ivanova (see \cite{I}), we say that a subset $S$ of a group $G$ is \textit{finitely $p$-separable} if for every $g$ $\in$ $G\setminus S$, there exists a homomorphism $\varphi$ from $G$ onto a finite $p$-group $P$ such that $\varphi(g)$ $\notin$ $\varphi(S)$. Note that $G$ is conjugacy $p$-separable if and only if every conjugacy class of $G$ is finitely $p$-separable. 

Examples of groups which are known to be conjugacy $p$-separable include free groups (see, for example, \cite{LS}) and fundamental groups of oriented closed surfaces (see \cite{P}).

There is a connection between these notions and a topology on $G$, the ``pro-$p$ topology" on $G$. The \textit{pro-$p$ topology} on $G$ is defined by taking the normal subgroups of $p$-power index in $G$ as a basis of neighbourhoods of 1 (see \cite{RZ}). Equipped with the pro-$p$ topology, $G$ becomes a topological group. Observe that $G$ is Hausdorff if and only if it is residually $p$-finite. One can show that a subset $S$ of $G$ is closed in the pro-$p$ topology on $G$ if and only if it is finitely $p$-separable. Thus, $G$ is conjugacy $p$-separable if and only if every conjugacy class of $G$ is closed in the pro-$p$ topology on $G$.
\\

In \cite{CZ}, Chagas and Zalesskii constructed an example of a conjugacy separable group possessing a non conjugacy separable subgroup of finite index. This led them to introduce the notion of ``hereditarily conjugacy separable group". A group $G$ is said to be \textit{hereditarily conjugacy separable} if every subgroup of finite index in $G$ is conjugacy separable. 

Recall that a \textit{subnormal subgroup} of a group $G$ is a subgroup $H$ of $G$ such that there exists a finite sequence of subgroups of $G$:
\begin{center}
$H$ = $H_0$ $<$ $H_1$ $<$ ... $<$ $H_n$ = $G$,
\end{center}
such that $H_i$ is normal in $H_{i+1}$ for all $i$ $\in$ \{0,...,$n-1$\}.

A subgroup $H$ of a group $G$ is open in the pro-$p$ topology on $G$ if and only if it is subnormal of $p$-power index (see Lemma \ref{A.1}). This leads us to the following definition, which naturally generalizes that of \cite{CZ}:

\begin{definition}
Let $G$ be a group. We say that $G$ is hereditarily conjugacy $p$-separable if every subnormal subgroup of $p$-power index in $G$ is conjugacy $p$-separable.
\end{definition}

In \cite{M}, Minasyan proved that right-angled Artin groups are hereditarily conjugacy separable. Our main theorem is the following:

\begin{reptheorem}{6.15}
Every right-angled Artin group is hereditarily conjugacy $p$-separable.
\end{reptheorem}

We will now discuss some applications of Theorem \ref{6.15}. The first application that we mention is an application of Theorem \ref{6.15} to separability properties of $G_\Gamma$:

\begin{repcorollary}{7.1}
Every right-angled Artin group is conjugacy separable in the class of torsion-free nilpotent groups.
\end{repcorollary}

Let $\mathcal{P}$ be a group property. A group $G$ is said to be \textit{virtually $\mathcal{P}$} if there exists a finite index subgroup $H$ $<$ $G$ such that $H$ has Property $\mathcal{P}$. Combining Theorem \ref{6.15} with a result of Paris (see \cite{P}), we obtain the following:

\begin{repcorollary}{7.4}
The outer automorphism group of a right-angled Artin group is virtually residually $p$-finite.
\end{repcorollary}

On the other hand, combining Theorem \ref{6.15} with a result of Myasnikov (see \cite{My}), we obtain the following:

\begin{repcorollary}{7.6}
The outer automorphism group of a right-angled Artin group is residually $\mathcal{K}$, where $\mathcal{K}$ is the class of all outer automorphism groups of finite $p$-groups.
\end{repcorollary}

The next application was suggested to the author by Ruth Charney and Luis Paris.

The natural action $Aut(G_\Gamma)$ $\rightarrow$ $GL_r(\mathbbm{Z})$ of $Aut(G_\Gamma)$ on $H_1(G_\Gamma,\mathbbm{Z})$ gives rise to a homomorphism $Out(G_\Gamma)$ $\rightarrow$ $GL_r(\mathbbm{Z})$, whose kernel is called the \textit{Torelli group} of $G_\Gamma$ -- by analogy with the Torelli group of a mapping class group. In Section 7, we combine well-known results of Bass-Lubotzky (see \cite{BL}), and Duchamp-Krob (see \cite{DK1}, \cite{DK2}) with Theorem \ref{6.15} to attain the following:

\begin{reptheorem}{7.14}
The Torelli group of a right-angled Artin group is residually torsion-free nilpotent.
\end{reptheorem}

\begin{repcorollary}{7.15}
The Torelli group of a right-angled Artin group is residually $p$-finite.
\end{repcorollary}

Recall that a group $G$ is said to be \textit{bi-orderable} if it can be endowed with a total order $\leq$ such that if $g$ $\leq$ $h$, then $kg$ $\leq$ $kh$ and $gk$ $\leq$ $hk$ for all $g$, $h$, $k$ $\in$ $G$. 

\begin{repcorollary}{7.16}
The Torelli group of a right-angled Artin group is bi-orderable.
\end{repcorollary}

Our proof follows closely that of Minasyan (see \cite{M}). Both proofs proceed by induction on the rank of $G_\Gamma$. The key observation is that a right-angled Artin group of rank $r$ can be written as an HNN extension of any of its special subgroups of rank $r-1$. After passing to an HNN extension of a finite group (which is known to be virtually free), Minasyan applies a theorem of Dyer stating that virtually free groups are conjugacy separable (see \cite{Dy1}).
\\

This paper is organized as follows. In Section 3, we introduce the $p$-centralizer condition which is the analogue of the centralizer condition in \cite{M}, and we prove that a group is hereditarily conjugacy $p$-separable if and only if it is conjugacy $p$-separable and satisfies the $p$-centralizer condition.  In Section 4, we prove the following analogue of Dyer's theorem for conjugacy $p$-separability:

\begin{theorem}
Every extension of a free group by a finite $p$-group is conjugacy $p$-separable.
\end{theorem}

Section 5 deals with retractions that are key tools in the proof of our main theorem, which is the object of Section 6.
\\

My gratefulness goes to my Ph.D. thesis advisor, Luis Paris, for his trust, time and advice. I am in debt to Ashot Minasyan for pointing out a mistake in an earlier draft of this work, and for directing me to the paper of Aschenbrenner and Friedl \cite{AF2}. I also wish to thank the referee for his (or her) many useful suggestions.

\section{HNN extensions}

\hspace{5mm}In this section, we recall the definition and basic properties of HNN extensions (see \cite{LS}).
\\

Let $H$ be a group. Then by the notation:
\begin{center}
$\langle$ $H$, $s$,... $\mid$ $r$,... $\rangle$,
\end{center}
we mean the group defined by the presentation whose generators are the generators of $H$ together with $s$,... and the relators of $H$ together with $r$,...

Let $H$ be a group, and let $K$ be a subgroup of $H$. The \textit{HNN extension of $H$ relative to $K$} is the group defined by the presentation:

\begin{center}
$G$ = $\langle$ $H$, $t$ $\mid$ $t^{-1}kt$ = $k$, $\forall$ $k$ $\in$ $K$ $\rangle$.
\end{center}

Every element of $G$ can be written as a word $x_0t^{a_1}x_1\cdots t^{a_n}x_n$ ($n$ $\geq$ 0, $x_0$,...,$x_n$ $\in$ $H$, $a_1$,...,$a_n$ $\in$ $\mathbbm{Z}\setminus\{0\}$). Following Minasyan (see \cite{M}), we will say that the word $x_0t^{a_1}x_1\cdots t^{a_n}x_n$ is \textit{reduced} if $x_0$ $\in$ $H$, $x_1$,...,$x_{n-1}$ $\in$ $H\setminus K$, and $x_n$ $\in$ $H$. Every element of $G$ can be written as a reduced word. Note that our definition of a reduced word is stronger than the definition of a reduced word in \cite{LS}.

\begin{lemma}[Britton's Lemma]
If a word $x_0t^{a_1}x_1\cdots t^{a_n}x_n$ is reduced with $n$ $\geq$ 1, then $x_0t^{a_1}x_1\cdots t^{a_n}x_n$ $\neq$ 1.
\end{lemma}

\textit{Proof}: Proved in \cite{LS} (see Theorem IV.2.1).\hfill$\square$

\begin{lemma}
If $x_0t^{a_1}x_1\cdots t^{a_n}x_n$ and $y_0t^{b_1}y_1\cdots t^{b_m}y_m$ are reduced words such that $x_0t^{a_1}x_1\cdots t^{a_n}x_n$ = $y_0t^{b_1}y_1\cdots t^{b_m}y_m$, then $m$ = $n$ and $a_i$ = $b_i$ for all $i$ $\in$ \{1,...,$n$\}.
\end{lemma}

\textit{Proof}: Proved in \cite{LS} (see Lemma IV.2.3).\hfill$\square$
\\

A \textit{cyclic permutation} of the word $t^{a_1}x_1\cdots t^{a_n}x_n$ is a word $t^{a_k}x_k\cdots t^{a_n}x_n\\ t^{a_1}x_1\cdots t^{a_{k-1}}x_{k-1}$ with $k$ $\in$ \{1,...,$n$\}. A word $t^{a_1}x_1\cdots t^{a_n}x_n$ is said to be \textit{cyclically reduced} if any cyclic permutation of $t^{a_1}x_1\cdots t^{a_n}x_n$ is reduced. Note that, if $t^{a_1}x_1\cdots t^{a_n}x_n$ is reduced and $n$ $\geq$ 2, then $t^{a_1}x_1\cdots t^{a_n}x_n$ is cyclically reduced if and only if $x_n$ $\in$ $H\setminus K$. Every element of $G$ is conjugate to a cyclically reduced word.

\begin{lemma}[Collins' Lemma] \label{2.3}
If $g$ = $t^{a_1}x_1\cdots t^{a_n}x_n$ ($n$ $\geq$ 1) and $h$ = $t^{b_1}y_1\cdots t^{b_m}y_m$ ($m$ $\geq$ 1) are cyclically reduced and conjugate, then there exists a cyclic permutation $h^*$ of $h$ and an element $\alpha$ $\in$ $K$ such that $g$ = $\alpha h^*\alpha^{-1}$.
\end{lemma}

\textit{Proof}: Proved in \cite{LS} (see Theorem IV.2.5).\hfill$\square$
\\

\begin{remark}
There exists a natural homomorphism $f$ : $G$ $\rightarrow$ $H$, defined by $f(h)$ = $h$ for all $h$ $\in$ $H$, and $f(t)$ = 1.
\end{remark}

\begin{remark}\label{2.5}
Let $P$ be a group and let $\varphi$ : $H$ $\rightarrow$ $P$ be a homomorphism. Let $Q$ be the HNN extension of $P$ relative to $\varphi(K)$:
\begin{center}
$Q$ = $\langle$ $P$, $\overline{t}$ $\mid$ $\overline{t}^{-1}\varphi(k)\overline{t}$ = $\varphi(k)$, $\forall$ $k$ $\in$ $K$ $\rangle$.
\end{center}
Then $\varphi$ induces a homomorphism $\overline{\varphi}$ : $G$ $\rightarrow$ $Q$, defined by $\overline{\varphi}(h)$ = $\varphi(h)$ for all $h$ $\in$ $H$, and $\overline{\varphi}(t)$ = $\overline{t}$.
\end{remark}

\begin{lemma}
With the notations of Remark \ref{2.5}, $ker(\overline{\varphi})$ is the normal closure of $ker(\varphi)$ in $G$.
\end{lemma}

\textit{Proof}: Proved in \cite{M} (see Lemma 7.5).\hfill$\square$
\\

The following observation is the key in the proof of our main theorem:

\begin{remark}
Let $G$ be a right-angled Artin group of rank $r$ ($r$ $\geq$ 1). Let $H$ be a special subgroup of $G$ of rank $r-1$. In other words, there is a partition of $V$, $V$ = $W$$\cup$\{$t$\}, such that $H$ = $\langle W \rangle$. Then $G$ can be written as the HNN extension of $H$ relative to the special subgroup $K$ = $\langle link(t) \rangle$ of $H$:
\begin{center}
$G$ = $\langle$ $H$, $t$ $\mid$ $t^{-1}kt$ = $k$, $\forall$ $k$ $\in$ $K$ $\rangle$.
\end{center}
\end{remark}

\section{Hereditary conjugacy $p$-separability and $p$-cen\-tralizer condition}

\hspace{5mm}We start with an observation that the reader has to keep in mind, because it will be used repeatedly in the rest of the paper: if $H$ and $K$ are two normal subgroups of $p$-power index in a group $G$, then $H\cap K$ is a normal subgroup of $p$-power index in $G$.
\\

The centralizer conditon was first introduced by Chagas and Zalesskii as a sufficient condition for a conjugacy separable group to be hereditarily conjugacy separable (see \cite{CZ}). Thereafter Minasyan showed that this condition is also necessary; that is, a group is hereditarily conjugacy separable if and only if it is conjugacy separable and satisfies the centralizer condition (see \cite{M}). We make the following definition, which naturally generalizes that of \cite{M}: 

\begin{definition}
We say that $G$ satisfies the $p$-centralizer condition ($pCC$) if, for every normal subgroup $H$ of $p$-power index in $G$, and for all $g$ $\in$ $G$, there exists a normal subgroup $K$ of $p$-power index in $G$ such that $K$ $<$ $H$, and:
\begin{center}
$C_{G/K}(\varphi(g))$ $\subset$ $\varphi(C_G(g)H)$,
\end{center}
where $\varphi$ : $G$ $\rightarrow$ $G/K$ denotes the canonical projection.
\end{definition}

We shall show that a group $G$ is hereditarily conjugacy $p$-separable if and only if it is conjugacy $p$-separable and satisfies the $p$-centralizer condition (see Proposition \ref{3.6}). If $H$ is a subgroup of $G$, and $g$ $\in$ $G$, we set $C_H(g)$ = \{$h$ $\in$ $H$ $\mid$ $gh$ = $hg$\}. For technical reasons, we have to introduce the following
definitions:

\begin{definition}
Let $G$ be a group, $H$ be a subgroup of $G$, and $g$ $\in$ $G$. We say that the pair $(H,g)$ satisfies the $p$-centralizer condition in $G$ ($pCC_G$) if, for every normal subgroup $K$ of $p$-power index in $G$, there exists a normal subgroup $L$ of $p$-power index in $G$ such that $L$ $<$ $K$, and:
\begin{center}
$C_{\varphi(H)}(\varphi(g))$ $\subset$ $\varphi(C_H(g)K)$,
\end{center}
where $\varphi$ : $G$ $\rightarrow$ $G/L$ denotes the canonical projection. We say that $H$ satisfies the $p$-centralizer condition in $G$ ($pCC_G$) if the pair $(H,g)$ satisfies the $p$-centralizer condition in $G$ for all $g$ $\in$ $G$.
\end{definition}

If $G$ is a group, $H$ is a subgroup of $G$, and $g$ $\in$ $G$, we set: $g^H$ = \{$\alpha$$g\alpha^{-1}$ $\mid$ $\alpha$ $\in$ $H$\}. In order to prove Proposition \ref{3.6}, we need the following statements, which are the analogues of some statements obtained in \cite{M} (Lemma 3.4, Corollary 3.5, and Lemma 3.7 respectively):

\begin{lemma} \label{3.3}
Let $G$ be a group, $H$ be a subgroup of $G$, and $g$ $\in$ $G$. Suppose that the pair $(G,g)$  satisfies $pCC_G$, and that $g^G$ is finitely $p$-separable in $G$. If $C_G(g)H$ is finitely $p$-separable in $G$, then $g^H$ is also finitely $p$-separable in $G$.
\end{lemma}

\textit{Proof}: Let $h$ $\in$ $G$ such that $h$ $\notin$ $g^H$. If $h$ $\notin$ $g^G$, then, since $g^G$ is finitely $p$-separable in $G$, there exists a homomorphism $\varphi$ from $G$ onto a finite $p$-group $P$ such that $\varphi(h)$ $\notin$ $\varphi(g^G)$. In particular, $\varphi(h)$ $\notin$ $\varphi(g^H)$. Thus we can assume that $h$ $\in$ $g^G$. Let $\alpha$ $\in$ $G$ be such that $h$ = $\alpha$$g\alpha^{-1}$. Suppose that $C_G(g)\cap$$\alpha^{-1}H$ $\neq$ $\emptyset$. Let $k$ $\in$ $C_G(g)\cap$$\alpha^{-1}H$. Then $\alpha$$k$ $\in$ $H$, and $h$ = $\alpha$$g\alpha^{-1}$ = $\alpha$$kg(\alpha$$k)^{-1}$ $\in$ $g^H$ -- which is a contradiction. Thus $C_G(g)\cap$$\alpha^{-1}H$ = $\emptyset$, i.e. $\alpha^{-1}$ $\notin$ $C_G(g)H$. As $C_G(g)H$ is finitely $p$-separable in $G$, there exists a normal subgroup $K$ of $p$-power index in $G$ such that $\alpha^{-1}$ $\notin$ $C_G(g)HK$. Now the condition $pCC_G$ implies that there exists a normal subgroup $L$ of $p$-power index in $G$ such that $L$ $<$ $K$, and:
\begin{center}
$C_{G/L}(\varphi(g))$ $\subset$ $\varphi(C_G(g)K)$,
\end{center}
where $\varphi$ : $G$ $\rightarrow$ $G/L$ denotes the canonical projection. We claim that $\varphi(h)$ $\notin$ $\varphi(g^H)$. Indeed, if there is $\beta$ $\in$ $H$ such that $\varphi(h)$ = $\varphi$($\beta$$g$$\beta^{-1}$), then $\varphi(\alpha^{-1}\beta)\varphi(g)$ = $\varphi(\alpha^{-1}\beta)\varphi(\beta^{-1}h\beta)$ = $\varphi(\alpha^{-1}h\alpha)\varphi(\alpha^{-1}\beta)$ = $\varphi(g)\varphi(\alpha^{-1}\beta)$, i.e. $\varphi(\alpha^{-1}\beta)$ $\in$ $C_{G/L}(\varphi(g))$. But then $\varphi(\alpha^{-1})$ $\in$ $C_{G/L}(\varphi(g))\varphi(H)$ $\subset$ $\varphi(C_G(g)\\KH)$. Hence $\alpha^{-1}$ $\in$ $C_G(g)HKL$ = $C_G(g)HK$ (because $L$ $<$ $K$) -- which is a contradiction.
\hfill$\square$

\begin{corollary} \label{3.4}
Let $G$ be a conjugacy $p$-separable group satisfying $pCC$, and $H$ be a subgroup of $G$ such that $C_G(h)H$ is finitely $p$-separable in $G$ for all $h$ $\in$ $H$. Then $H$ is conjugacy $p$-separable. Moreover, for all $h$ $\in$ $H$, $h^H$ is finitely $p$-separable in $G$.
\end{corollary}

\textit{Proof}: Let $h$ $\in$ $H$. Since $G$ satisfies $pCC$, the pair $(G,h)$ satisfies $pCC_G$. Since $G$ is conjugacy $p$-separable, $h^G$ is finitely $p$-separable in $G$. Lemma \ref{3.3} now implies that $h^H$ is finitely $p$-separable in $G$. Therefore $h^H$ is finitely $p$-separable in H.\hfill$\square$

\begin{lemma} \label{3.5}
Let $G$ be a group, $H$ be a subgroup of $G$, and $g$ $\in$ $G$. Let $K$ be a normal subgroup of $p$-power index in $G$. If $g^{H \cap K}$ is finitely $p$-separable in $G$, then there exists a normal subgroup $L$ of $p$-power index in $G$ such that $L$ $<$ $K$, and:
\begin{center}
$C_{\varphi(H)}(\varphi(g))$ $\subset$ $\varphi(C_H(g)K)$,
\end{center}
where $\varphi$ : $G$ $\rightarrow$ $G/L$ denotes the canonical projection.
\end{lemma}

\textit{Proof}: Note that $H\cap$$K$ is of finite index $n$ in $H$. Actually, $H\cap$$K$ is of $p$-power index in $H$ (because $\frac{H}{H\cap K}$ $\simeq$ $\frac{KH}{K}$ $<$ $\frac{G}{K}$), but this is not needed here. There exist $\alpha_1$,...,$\alpha_n$ $\in$ $H$ such that $H$ = $\sqcup_{i=1}^n\alpha_i(H\cap K)$. Up to renumbering, we can assume that there exists $l$ $\in$ \{0,...,$n$\} such that $\alpha_i^{-1}g\alpha_i$ $\in$ $g^{H\cap K}$ for all $i$ $\in$ \{1,...,$l$\} and $\alpha_i^{-1}g\alpha_i$ $\notin$ $g^{H\cap K}$ for all $i$ $\in$ \{$l$+1,...,$n$\}. By the assumptions, there exists a normal subgroup $L$ of $p$-power index in $G$ such that $\alpha_i^{-1}g\alpha_i$ $\notin$ $g^{H\cap K}L$ for all $i$ $\in$ \{$l$+1,...,$n$\}. Up to replacing $L$ by $L\cap$$K$, we can assume that $L$ $<$ $K$. Let $\varphi$ : $G$ $\rightarrow$ $G/L$ be the canonical projection. Let $\overline{h}$ $\in$ $C_{\varphi(H)}(\varphi(g))$. There exists $h$ $\in$ $H$ such that $\overline{h}$ = $\varphi(h)$. There exist $i$ $\in$ \{1,...,$n$\} and $k$ $\in$ $H\cap$$K$ such that $h$ = $\alpha_ik$. We have $\varphi(h^{-1}gh)$ = $\varphi(h)^{-1}\varphi(g)\varphi(h)$ = $\varphi(g)$. Thus $h^{-1}gh$ $\in$ $gL$. But then $\alpha_i^{-1}g\alpha_i$ = $kh^{-1}ghk^{-1}$ $\in$ $kgLk^{-1}$ = $kgk^{-1}L$ $\subset$ $g^{H\cap K}L$. Therefore $i$ $\leq$ $l$. Then there exists $\beta$ $\in$ $H\cap$$K$ such that: $\alpha_i^{-1}g\alpha_i$ = $\beta$$g\beta^{-1}$. This is to say that $\alpha_i\beta$ $\in$ $C_H(g)$, and then $h$ = $\alpha_ik$ = $(\alpha_i\beta)(\beta^{-1}k)$ $\in$ $C_H(g)(H\cap K)$ $\subset$ $C_H(g)K$. We have shown that $C_{\varphi(H)}(\varphi(g))$ $\subset$ $\varphi(C_H(g)K)$.\hfill$\square$
\\

We are now ready to prove:

\begin{proposition} \label{3.6}
A group is hereditarily conjugacy $p$-separable if and only if it is conjugacy $p$-separable and satisfies $pCC$.
\end{proposition}

\textit{Proof}: Suppose that $G$ is conjugacy $p$-separable and satisfies $pCC$. Let $H$ be a subnormal subgroup of $p$-power index in $G$. Thus $H$ is closed in the pro-$p$ topology on $G$ (because $G\setminus H$ = $\cup\{gH \mid g \notin H\}$). Let $h$ $\in$ $H$. The set $C_G(h)H$ is a finite union of left cosets modulo $H$ and thus is closed in the pro-$p$ topology on $G$. Corollary \ref{3.4} now implies that $H$ is conjugacy $p$-separable. Therefore $G$ is hereditarily conjugacy $p$-separable. Suppose now that $G$ is hereditarily conjugacy $p$-separable. In particular, $G$ is conjugacy $p$-separable. We shall show that $G$ satisfies $pCC$. Let $g$ $\in$ $G$. Let $K$ be a normal subgroup of $p$-power index in $G$. Let $H$ = $K \langle g \rangle$. Since $K$ $<$ $H$, $[G:H]$ is a power of $p$. As $\frac{G}{K}$ is a finite $p$-group, every subgroup of it is subnormal. Thus $H$ is subnormal in $G$. Therefore $H$ is conjugacy $p$-separable. Note that $g^{G\cap K}$ = $g^K$ = $g^H$ $\subset$ $H$. As $g^H$ is closed in the pro-$p$ topology on $H$, it is closed in the pro-$p$ topology on $G$, because the topology induced on $H$ by the pro-$p$ topology on $G$ coincides with the pro-$p$ topology on $H$ (see, for example, \cite{RZ2}, Corollary 5.8). The result now follows from lemma \ref{3.5}.\hfill$\square$

\section{Extensions of free groups by finite $p$-groups are conjugacy $p$-separable}

\hspace{5mm}We start with an observation that the reader has to keep in mind because it will be used repeatedly in the proof of Theorem \ref{4.2}: if $\varphi$ : $G$ $\rightarrow$ $H$ is a homomorphism from a group $G$ to a group $H$, whose kernel is torsion-free, then the restriction of $\varphi$ to any finite subgroup of $G$ is injective.
\\

We need the following lemma:

\begin{lemma}\label{4.1}
Let $G$ = $G_1\ast\cdots\ast G_n$ be a free product of $n$ conjugacy $p$-separable groups $G_1$,...,$G_n$. Let $g$, $h$ $\in$ $G\setminus\{1\}$ be two non-trivial elements of finite order in $G$ such that $g$ $\nsim$ $h$. There exists a homomorphism $\varphi$ from $G$ onto a finite $p$-group $P$ such that $\varphi(g)$ $\nsim$ $\varphi(h)$.
\end{lemma}

\textit{Proof}: Since $g$ is of finite order in $G$, there exists $i$ $\in$ \{1,...,$n$\} such that $g$ is conjugate to an element of finite order in $G_i$. Thus we may assume that $g$ belongs to $G_i$. Similarly, we may assume that there exists $j$ in \{1,...,$n$\} such that $h$ belongs to $G_j$. Suppose that $i$ $\neq$ $j$. Let $\varphi$ : $G_i$ $\rightarrow$ $P$ be a homomorphism from $G_i$ onto a finite $p$-group $P$ such that $\varphi(g)$ $\neq$ 1. Let $\widetilde{\varphi}$ : $G$ $\rightarrow$ $P$ be the natural homomorphism extending $\varphi$. Then $\widetilde{\varphi}(g)$ $\nsim$ $\widetilde{\varphi}(h)$. Suppose that $i$ = $j$. Then $g$ and $h$ are not conjugate in $G_i$ -- otherwise they would be conjugate in $G$. Since $G_i$ is conjugacy $p$-separable, there exists a homomorphism $\varphi$ : $G_i$ $\rightarrow$ $P$ from $G_i$ onto a finite $p$-group $P$ such that $\varphi(g)$ $\nsim$ $\varphi(h)$. Let $\widetilde{\varphi}$ : $G$ $\rightarrow$ $P$ be defined as above. We have $\widetilde{\varphi}(g)$ $\nsim$ $\widetilde{\varphi}(h)$.\hfill$\square$
\\

In Section 4, by a graph, we mean a unoriented graph, possibly with loops or multiple edges.
\\

Recall that a \textit{graph of groups} is a connected graph $\Gamma$ = $(V,E)$, together with a function $\mathcal{G}$ which assigns:
\begin{itemize}
\item to each vertex $v$ $\in$ $V$, a group $G_v$,
\item and to each edge $e$ = \{$v$,$w$\} $\in$ $E$, a group $G_e$ together with two injective homomorphisms $\alpha_e$ : $G_e$ $\rightarrow$ $G_v$ and $\beta_e$ : $G_e$ $\rightarrow$ $G_w$ -- we are not assuming that $v$ $\neq$ $w$ --,
\end{itemize}
(see \cite{Se}, see also \cite{Dy1}). The groups $G_v$ ($v$ $\in$ $V$) are called the \textit{vertex groups} of $\Gamma$, the groups $G_e$ ($e$ $\in$ $E$) are called the \textit{edge groups} of $\Gamma$. The monomorphisms $\alpha_e$ and $\beta_e$ ($e$ $\in$ $E$) are called the \textit{edge monomorphisms}. The images of the edge groups under the edge monomorphisms are called the \textit{edge subgroups}.

Choose disjoint presentations $G_v$ = $\langle$ $X_v$ $\mid$ $R_v$ $\rangle$ for the vertex groups of $\Gamma$. Choose a maximal tree $T$ in $\Gamma$. Assign a direction to each edge of $\Gamma$. Let \{$t_e$ $\mid$ $e$ $\in$ $E$\} be a set in one-to-one correspondence with the set of edges of $\Gamma$, and disjoint from the $X_v$ ($v$ $\in$ $V$). The \textit{fundamental group} of the above graph of groups $\Gamma$ is the group $G_\Gamma$ defined by the presentation whose generators are:
\begin{center}
$X_v$ ($v$ $\in$ $V$), \\
$t_e$ ($e$ $\in$ $E$)
\end{center}
(called vertex and edge generators respectively) and whose relations are:
\begin{center}
$R_v$ ($v$ $\in$ $V$), \\
$t_e$ = 1 ($e$ $\in$ $T$), \\
$t_e\alpha_e(g_e)t_e^{-1}$ = $\beta_e(g_e)$, $\forall$ $g_e$ $\in$
$G_e$ ($e$ $\in$ $E$).
\end{center}
(called vertex, tree, and edge relations respectively). One can prove that this is well-defined -- that is, independent of our choice of $T$, etc. Note that it suffices to write the edge relations for $g_e$ in a set of generators for $G_e$.
\\

\underline{Convention}: The groups $G_v$ ($v$ $\in$ $V$) and $G_e$ ($e$ $\in$ $E$) will be regarded as subgroups of $G_\Gamma$.
\\

Let $\{\Gamma_i\}_{i \in I}$ be a collection of connected and pairwise disjoint subgraphs of $\Gamma$. We may define a graph of groups $\Gamma^*$ from $\Gamma$ \textit{by contracting $\Gamma_i$ to a point for all $i$ $\in$ $I$}, as follows. The graph $\Gamma^*$ is obtained from $\Gamma$ by contracting $\Gamma_i$ to a point $p_i$ for all $i$ $\in$ $I$. The function $\mathcal{G^*}$ is obtained from $\mathcal{G}$ by using the fundamental group of $\Gamma_i$ for the vertex group at $p_i$, and by composing the edge monomorphisms of $\Gamma$ by the natural inclusions of the vertex groups of $\Gamma_i$ into the fundamental group of $\Gamma_i$, if necessary. Clearly, $G_\Gamma$ is isomorphic to the fundamental group $G_{\Gamma^*}$ of $\Gamma^*$.
\\

If $\pi$ : $G_\Gamma$ $\rightarrow$ $H$ is a homomorphism from $G_\Gamma$ to a group $H$, such that the restriction of $\pi$ to each edge subgroup of $\Gamma$ is injective, then we may define a graph of groups $\Gamma'$ from $\Gamma$ \textit{by composing with $\pi$}, as follows. The vertex set of $\Gamma'$ is $V$, and the edge set of $\Gamma'$ is $E$. The vertex groups of $\Gamma'$ are the groups $G_v'$ = $\pi(G_v)$ ($v$ $\in$ $V$), and the edge groups of $\Gamma'$ are the groups $G_e'$ = $G_e$ ($e$ $\in$ $E$). The edge monomorphisms are the monomorphisms $\alpha_e'$ = $\pi\circ\alpha_e$ and $\beta_e'$ = $\pi\circ\beta_e$ ($e$ $\in$ $E$). Present $G_\Gamma$ and $G_{\Gamma'}$ using the same symbols for edge generators and with the same choice of maximal tree. There exist two homomorphisms, $\pi_V$ : $G_\Gamma$ $\rightarrow$ $G_{\Gamma'}$ and $\pi_E$ : $G_{\Gamma'}$ $\rightarrow$ $H$, such that the diagram:
$$
\xymatrix{
  G_\Gamma \ar[d]_{\pi_V}  \ar[r]^\pi & H \\
  G_{\Gamma'} \ar[ru]_{\pi_E}
  }
$$
commutes, and that the restriction of $\pi_E$ to each vertex group of $G_{\Gamma'}$ is injective. The homomorphism $\pi_V$ is given by:
\begin{center}
$(\pi_V)_{\mid_{G_v}}$ = $\pi_{\mid_{G_v}}$, $\forall$ $v$ $\in$ $V$,\\
$\pi_V(t_e)$ = $t_e$, $\forall$ $e$ $\in$ $E$.
\end{center}
And the homomorphism $\pi_E$ is given by:
\begin{center}
$(\pi_E)_{\mid_{G_v'}}$ = $(id_H)_{\mid_{G_v'}}$, $\forall$ $v$ $\in$ $V$,\\
$\pi_E(t_e)$ = $\pi(t_e)$, $\forall$ $e$ $\in$ $E$.
\end{center}

\vspace{2mm}

In \cite{Dy1}, Dyer proved that every extension of a free group by a finite group is conjugacy separable. The following theorem is the analogue of Dyer's theorem for conjugacy $p$-separability.

\begin{theorem}\label{4.2}
Every extension of a free group by a finite $p$-group is conjugacy $p$-separable.
\end{theorem}

\textit{Proof}: Our proof is inspired by that of Dyer (see \cite{Dy1}). Let $G$ be an extension of a free group by a finite $p$-group. In other words, there exists a short exact sequence:
$$
\xymatrix{
    1 \ar[r] & F \ar[r] & G \ar[r]^\pi & P \ar[r] & 1
  },
$$
where $F$ is a free group, and $P$ is a finite $p$-group. Let $g$ $\in$ $G$. Let $h$ $\in$ $G$ such that $g$ $\nsim$ $h$.
\\

\underline{Step 1}: We show that we may assume that $G$ satisfies a short exact sequence:
$$
\xymatrix{
    1 \ar[r] & F \ar[r] & G \ar[r]^\pi & C_{p^n} \ar[r] & 1
  },
$$
where $F$ is a free group, $n$ $\geq$ 1, $C_{p^n}$ denotes the cyclic group of order $p^n$, and $\pi(g)$ = $\pi(h)$.
\\

Since $G$ is an extension of a free group by a finite $p$-group, $G$ is residually $p$-finite by \cite{G}, Lemma 1.5. Therefore, if $g$ = 1, then $g^G$ = \{1\} is finitely $p$-separable in $G$. On the other hand, if $g$ is of infinite order in $G$, then $g^G$ is finitely $p$-separable in $G$ by \cite{I}, Proposition 5. Therefore we may assume that $g$ $\neq$ 1 and that $g$ is of finite order in $G$. Similarly, we may assume that $h$ $\neq$ 1 and that $h$ is of finite order in $G$. If $\pi(g)$ and $\pi(h)$ are not conjugate in $P$, we are done. Thus, up to replacing $h$ by a conjugate of itself, we may assume that $\pi(g)$ = $\pi(h)$. Since $ker(\pi)$ = $F$ is torsion-free, $g$ and $h$ have the same order $p^n$ ($n$ $\in$ $\mathbbm{N}^*$). Let $H$ = $F\langle g\rangle$. Note that $H$ is a subgroup of $p$-power index in $G$, and that $g$ and $h$ belong to $H$. As $\frac{G}{F}$ = $P$ is a finite $p$-group, every subgroup of it is subnormal. Thus $H$ is subnormal in $G$. Then we may replace $G$ by $H$, by \cite{I}, Proposition 4\footnote{Strictly speaking, it follows from the proof of \cite{I}, Proposition 4, that, if there exists a homomorphism $\varphi$ : $H$ $\rightarrow$ $P$ from $H$ onto a finite $p$-group $P$ such that $\varphi(g)$ $\nsim$ $\varphi(h)$, then there exists a homomorphism $\psi$ : $G$ $\rightarrow$ $Q$ from $G$ onto a finite $p$-group $Q$ such that $\psi(g)$ $\nsim$ $\psi(h)$. The exact statement of \cite{I}, Proposition 4, is slightly different.}, so as to assume that $G$ satisfies the short exact sequence:
$$
\xymatrix{
    1 \ar[r] & F \ar[r] & G \ar[r]^\pi & C_{p^n} \ar[r] & 1
  }.
$$
\\

Now, $G$ is the fundamental group of a graph of groups $\Gamma$, whose vertex groups are all finite groups, by \cite{S}, Theorem. As $\pi_{\mid_{G_v}}$ is injective for all $v$ $\in$ $V$, $G_v$ is isomorphic to a subgroup of $C_{p^n}$ for all $v$ $\in$ $V$. From now on, the groups $G_v$ ($v$ $\in$ $V$) will be regarded as subgroups of $C_{p^n}$.
\\

\underline{Step 2}: We show that we may assume that all edge groups are non-trivial, that if two different vertices are connected by an edge, then the corresponding edge group is a proper subgroup of $C_{p^n}$, and that $g$ and $h$ belong to two different vertex groups.
\\

First, we show that we may assume that all edge groups are non-trivial. Indeed, Let $\Gamma_0$ be the subgraph of $\Gamma$ whose vertices are all the vertices of $\Gamma$, and whose edges are the edges of $\Gamma$ for which the edge group is non-trivial. Let $\Gamma_1$,...,$\Gamma_r$ be the connected components of $\Gamma_0$. Let $\Gamma^*$ be the graph of groups obtained from $\Gamma$ by contracting $\Gamma_i$ to a point for all $i$ $\in$ \{1,...,$r$\}. Let $T$ be a maximal tree of $\Gamma^*$. Then $G$ is isomorphic to the fundamental group $G^*$ of $\Gamma^*$. Observe that $G^*$ is the free product of the free group on \{$t_e$ $\mid$ $e$ $\in$ $E \setminus T$\} and the fundamental groups of the $\Gamma_i$ ($i$ $\in$ \{1,...,$r$\}). Thus, it suffices to consider the case where $\Gamma$ = $\Gamma_i$ ($i$ $\in$ \{1,...,$r$\}), by Lemma \ref{4.1}. Since each $\Gamma_i$ ($i$ $\in$ \{1,...,$r$\}) is a graph of groups whose edge groups are all non-trivial, the first part of the assertion is proved.

Now, we show that we may assume that if two different vertices are connected by an edge, then the corresponding edge group is a proper subgroup of $C_{p^n}$. Indeed, let $\Gamma_0$ be the subgraph of $\Gamma$ whose vertices are all the vertices of $\Gamma$, and whose edges are the edges of $\Gamma$ for which the edge group is isomorphic to $C_{p^n}$. Let $\Gamma_1$,...,$\Gamma_r$ be the connected components of $\Gamma_0$. Choose a maximal tree $T_i$ in $\Gamma_i$, for all $i$ $\in$ \{1,...,$r$\}. Let $\Gamma^*$ be the graph of groups obtained from $\Gamma$ by contracting $T_i$ to a point for all $i$ $\in$ \{1,...,$r$\}. Then $G$ is isomorphic to the fundamental group $G^*$ of $\Gamma^*$. Note that a vertex group of $\Gamma^*$ is either a vertex group of $\Gamma$, or the fundamental group of $T_i$, for some $i$ $\in$ \{1,...,$r$\}, in which case it is isomorphic to $C_{p^n}$ (because each $T_i$ ($i$ $\in$ \{1,...,$r$\}) is a tree of groups whose vertex and edge groups are all equal to $C_{p^n}$). Thus, we may replace $\Gamma$ by $\Gamma^*$, so that the second part of the assertion is proved.

Since $g$ is of finite order in $G$, there exists a vertex $v$ of $\Gamma$, an element $g_0$ of finite order in the vertex group $G_v$ of $v$, and an element $\alpha$ of $G$ such that $g$ = $\alpha g_0 \alpha^{-1}$. Similarly, there exists a vertex $w$ of $\Gamma$, an element $h_0$ of finite order in the vertex group $G_w$ of $w$, and an element $\beta$ of $G$ such that $h$ = $\beta h_0 \beta^{-1}$. As $C_{p^n}$ is abelian, we have: $\pi(g_0)$ = $\pi(h_0)$. Thus, up to replacing $g$ by $g_0$ and $h$ by $h_0$, we may assume that $g$ belongs to $G_v$, and $h$ belongs to $G_w$. Since $\pi_{\mid_{G_v}}$ is injective, and $\pi(g)$ = $\pi(h)$, we have $v$ $\neq$ $w$.
\\

\underline{Step 3}: We show that we may assume that $\Gamma$ has exactly two vertices, and that all edges join these two vertices.
\\

Indeed, choose a maximal tree $T$ in $\Gamma$. There is a path $P$ in $T$ joining $v$ to $w$. Choose an edge $e$ on this path. Then $T\setminus\{e\}$ is the disjoint union of two trees, $T_v$ and $T_w$ -- with $v$ $\in$ $T_v$ and $w$ $\in$ $T_w$. Let $\Gamma_v$ be the full subgraph of $\Gamma$ generated by the vertices of $T_v$, and $\Gamma_w$ be the full subgraph of $\Gamma$ generated by the vertices of $T_w$. Let $\Gamma^*$ be the graph of groups obtained from $\Gamma$ by contracting $\Gamma_v$ to a point $v^*$ and $\Gamma_w$ to a point $w^*$. Observe that $\Gamma^*$ has exactly two vertices and that all edges join these two vertices. The vertex groups of $\Gamma^*$ are the fundamental groups of $\Gamma_v$ and $\Gamma_w$, respectively. The edge groups of $\Gamma^*$ are non-trivial proper subgroups of $C_{p^n}$. And $G$ is isomorphic to the fundamental group $G^*$ of $\Gamma^*$. Now, since the restriction of $\pi$ to each edge subgroup of $\Gamma^*$ is injective, we may define a graph of groups $\Gamma'$ from $\Gamma^*$ by composing with $\pi$, as described above. Denote by $G'$ the fundamental group of $\Gamma'$. There exist two homomorphisms, $\pi_V$ : $G$ $\rightarrow$ $G'$ and $\pi_E$ : $G'$ $\rightarrow$ $C_{p^n}$, such that the diagram:
$$
\xymatrix{
  G \ar[d]_{\pi_V}  \ar[r]^\pi & C_{p^n} \\
  G' \ar[ru]_{\pi_E}
  }
$$
commutes, and that the restriction of $\pi_E$ to each vertex group of $\Gamma'$ is injective. Consequently, $ker(\pi_E)$ is free by \cite{Se}, II, 2.6., Lemma 8.

Set $g'$ = $\pi_V(g)$, and $h'$ = $\pi_V(h)$. As $g'$ and $h'$ have order $p^n$, the vertex groups of $\Gamma'$ are equal to $C_{p^n}$. The edge groups of $\Gamma'$ are non-trivial proper subgroups of $C_{p^n}$. Observe that $g'$ and $h'$ belong to two different vertex groups, and that $g'$ (resp. $h'$) is not conjugate to an element of one of the edge groups. Let $e$ be an edge of $\Gamma'$. Then $g'$ and $h'$ are not conjugate in $G_v'\ast_{G_e'}G_w'$, by \cite{MKS}, Theorem 4.6. Observe that $G'$ is an HNN extension (in the general sense) of $G_v'\ast_{G_e'}G_w'$ with stable letters $t_a$ ($a$ $\in$ $E\setminus\{e\}$), and associated subgroups $\alpha_a'(G_a')$ and $\beta_a'(G_a')$ ($a$ $\in$ $E\setminus\{e\}$). Therefore $g'$ and $h'$ are not conjugate in $G'$ (see, for example, \cite{Dy2}, Theorem 3). Thus, we may replace $\Gamma$ by $\Gamma'$, $G$ by $G'$, $g$ by $g'$, and $h$ by $h'$, so as to assume that $\Gamma$ has two vertices and that all edges join these two vertices.
\\

\underline{Step 4}: We show that we may assume that $\Gamma$ has at most two edges.
\\

Suppose that $\Gamma$ has more than two edges. Choose a maximal tree $T$ in $\Gamma$ -- that is, an edge of $\Gamma$. Present $G_v$ = $\langle g \mid g^{p^n} = 1 \rangle$, $G_w$ = $\langle h \mid h^{p^n} = 1 \rangle$, and $G$ as described above. Choose an edge $e$ $\in$ $E \setminus T$.
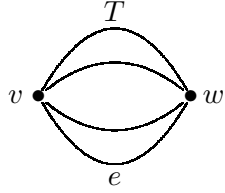
\[
\xygraph{ 
!{<0cm,0cm>;<1cm,0cm>:<0cm,1cm>::}
!{(0,0)}*{\bullet}="v"
!{(-0.3,0)}*\txt{$v$}
!{(2,0)}*{\bullet}="w"
!{(2.3,0)}*\txt{$w$}
"v"-@/^0.9cm/"w"^{\txt{$T$}} "v"-@/^0.45cm/"w" 
"v"-@/_0.45cm/"w" "v"-@/_0.9cm/"w"_{\txt{$e$}}
}
\]
The edge relations corresponding to $e$ can be reduced to the following:
\begin{center}
$t_e\alpha_e(g_e)t_e^{-1}$ = $\beta_e(g_e)$,
\end{center}
where $g_e$ is a generator of $G_e$. Let $p^s$ be the order of $G_e$ ($s$ $\in$ \{1,...,$n-1$\}). Then $\alpha_e(g_e)$ generates a subgroup of order $p^s$ of $G_v$. But there exists a unique subgroup of order $p^s$ in $G_v$; it is cyclic, generated by $g^{p^r}$, where $r$ = $n-s$. Thus, up to replacing $g_e$ by the preimage of $g^{p^r}$ under $\alpha_e$, we may assume that $\alpha_e(g_e)$ = $g^{p^r}$. There exists $k$ $\in$ $\mathbbm{N}$, such that $p$ and $k$ are coprime, and that $\beta_e(g_e)$ = $h^{kp^r}$. Therefore the edge relation corresponding to $e$ can be written:
\begin{center}
$t_eg^{p^r}t_e^{-1}$ = $h^{kp^r}$,
\end{center}
where $r$ $\in$ \{1,...,$n-1$\}, $k$ $\in$ $\mathbbm{N}$, and $p$ and $k$ are coprime. Now, since $\pi$ : $G$ $\rightarrow$ $C_{p^n}$ satisfies $\pi(g)$ = $\pi(h)$, we have: $\pi(g)^{p^r}$ = $\pi(h)^{kp^r}$ = $\pi(g)^{kp^r}$, and then $\pi(g)^{(k-1)p^r}$ = 1 (in $C_{p^n}$). As $\pi(g)$ has order $p^n$ in $C_{p^n}$, we deduce that $p^{n-r}$ divides $k - 1$. There exists $a$ $\in$ $\mathbbm{Z}$ such that $k$ = $ap^{n-r}$ + 1. We conclude that the edge relation corresponding to $e$ can be written:
\begin{center}
$t_eg^{p^r}t_e^{-1}$ = $h^{p^r}$,
\end{center}
where $r$ $\in$ \{1,...,$n-1$\}.
\\

Let $H$ be the normal subgroup of $G$ generated by the elements:
\begin{center}
$g$, $h$, $t_a$ ($a$ $\in$ $E\setminus\{e\}$), $t_e^p$.
\end{center}
Then $H$ has index p in $G$, and $g$ and $h$ belong to $H$. Thus we may replace $G$ by $H$ by \cite{I}, Proposition 4. Let $G_0$ be the fundamental group of the graph of groups $\Gamma\setminus\{e\}$. Set $G_0$ = $\langle X_0 \mid R_0 \rangle$, where the presentation is as fundamental group of the graph of groups $\Gamma\setminus\{e\}$. Set $G_i$ = $t_e^iG_0t_e^{-i}$ = $\langle X_i \mid R_i \rangle$, for all $i$ $\in$ \{1,...,$p-1$\}. Clearly \{1,$t_e$,...,$t_e^{p-1}$\} is a Schreier transversal for $H$ in $G$. The Reidemeister-Schreier method yields the presentation:

\begin{center}
$H$ = $\langle$ $X_0$, $X_1$,..., $X_{p-1}$, $u$ $\mid$ $R_0$, $R_1$,..., $R_{p-1}$, $g_1^{p^r}$ = $h_0^{p^r}$, $g_2^{p^r}$ = $h_1^{p^r}$,..., $g_{p-1}^{p^r}$ = $h_{p-2}^{p^r}$, $ug_0^{p^r}u^{-1}$ = $h_{p-1}^{p^r}$ $\rangle$,
\end{center}
where $u$ = $t_e^p$, $g_i$ = $t_e^igt_e^{-i}$ ($i$ $\in$ \{0,...,$p-1$\}), and $h_j$ = $t_e^jht_e^{-j}$ ($j$ $\in$ \{0,...,$p-1$\}). Replace $g$ by $g_0$, and $h$ by $h_1$. Observe that $H$ is the fundamental group of a graph of groups $\widetilde{\Gamma}$, as follows. The graph $\widetilde{\Gamma}$ has 2$p$ vertices, say $v_0$, $w_0$, $v_1$, $w_1$,...,$v_{p-1}$, $w_{p-1}$, and $p\left|E\right|$ edges. Let $\widetilde{\Gamma_i}$ be the full subgraph of $\widetilde{\Gamma}$ generated by \{$v_i$, $w_i$\} for all $i$ $\in$ \{0,...,$p-1$\}. Then $\widetilde{\Gamma_i}$ is isomorphic to $\Gamma\setminus\{e\}$. There is one edge joining $w_0$ to $v_1$, one edge joining $w_1$ to $v_2$,..., one edge joining $w_{p-2}$ to $v_{p-1}$, and one edge joining $v_0$ to $w_{p-1}$, and the edge groups associated to these egdes are isomorphic to $G_e$. Note that $g$ belongs to the vertex group of $v_0$ and $h$ belongs to the vertex group of $w_1$.
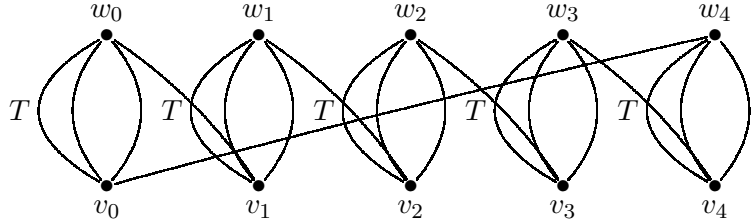
\[
\xygraph{ 
!{<0cm,0cm>;<1cm,0cm>:<0cm,1cm>::}
!{(-19,-1)}*{\bullet}="v_0"
!{(-19,-1.3)}*\txt{$v_0$}
!{(-19,1)}*{\bullet}="w_0"
!{(-19,1.3)}*\txt{$w_0$}
!{(-17,-1)}*{\bullet}="v_1"
!{(-17,-1.3)}*\txt{$v_1$}
!{(-17,1)}*{\bullet}="w_1"
!{(-17,1.3)}*\txt{$w_1$}
!{(-15,-1)}*{\bullet}="v_2"
!{(-15,-1.3)}*\txt{$v_2$}
!{(-15,1)}*{\bullet}="w_2"
!{(-15,1.3)}*\txt{$w_2$}
!{(-13,-1)}*{\bullet}="v_3"
!{(-13,-1.3)}*\txt{$v_3$}
!{(-13,1) }*{\bullet}="w_3"
!{(-13,1.3)}*\txt{$w_3$}
!{(-11,-1)}*{\bullet}="v_4"
!{(-11,-1.3)}*\txt{$v_4$}
!{(-11,1)}*{\bullet}="w_4"
!{(-11,1.3)}*\txt{$w_4$}
"v_0"-@/^0.9cm/"w_0"^{\txt{$T$}} "v_0"-@/^0.45cm/"w_0" "v_0"-@/_0.45cm/"w_0"
"v_1"-@/^0.9cm/"w_1"^{\txt{$T$}} "v_1"-@/^0.45cm/"w_1" "v_1"-@/_0.45cm/"w_1"
"v_2"-@/^0.45cm/"w_2" "v_2"-@/_0.45cm/"w_2" "v_2"-@/^0.9cm/"w_2"^{\txt{$T$}}
"v_3"-@/^0.45cm/"w_3" "v_3"-@/_0.45cm/"w_3" "v_3"-@/^0.9cm/"w_3"^{\txt{$T$}}
"v_4"-@/^0.9cm/"w_4"^{\txt{$T$}} "v_4"-@/^0.45cm/"w_4" "v_4"-@/_0.45cm/"w_4"
"w_0"-@/^0.15cm/"v_1" "w_1"-@/^0.15cm/"v_2" "w_2"-@/^0.15cm/"v_3"
"w_3"-@/^0.15cm/"v_4" "v_0"-@/^0.05cm/"w_4"
}
\]
Let $\Gamma^*$ be the graph of groups obtained from $\widetilde{\Gamma}$ by contracting $\widetilde{\Gamma_i}$ to a point for all $i$ $\in$ \{1,...,$p-1$\}. Then $G$ is isomorphic to the fundamental group of $\Gamma^*$. The graph $\Gamma^*$ has $p$ vertices, say $v_0$,...,$v_{p-1}$. There is one edge joining $v_0$ to $v_1$, one edge joining $v_1$ to $v_2$,..., one edge joining $v_{p-2}$ to $v_{p-1}$, and one edge joining $v_0$ to $v_{p-1}$, and the edge groups associated to these edges are all isomorphic to $G_e$. Note that $g$ belongs to the vertex group of $v_0$ and $h$ belongs to the vertex group of $v_1$.
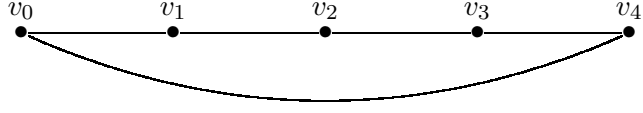
\[
\xygraph{
!{<0cm,0cm>;<1cm,0cm>:<0cm,1cm>::}
!{(-19,0)}*{\bullet}="v_0"
!{(-19,0.3)}*\txt{$v_0$}
!{(-17,0)}*{\bullet}="v_1"
!{(-17,0.3)}*\txt{$v_1$}
!{(-15,0)}*{\bullet}="v_2"
!{(-15,0.3)}*\txt{$v_2$}
!{(-13,0)}*{\bullet}="v_3"
!{(-13,0.3)}*\txt{$v_3$}
!{(-11,0)}*{\bullet}="v_4"
!{(-11,0.3)}*\txt{$v_4$}
"v_0"-"v_1" "v_1"-"v_2"
"v_2"-"v_3" "v_3"-"v_4" 
"v_0"-@/_0.9cm/"v_4"
}
\]
Let $T$ be the maximal tree $T$ = $v_0v_1\cdots v_{p-2}v_{p-1}$. Then $T\setminus\{v_0v_1\}$ is the disjoint union of two trees : $v_0$ and $v_1v_2\cdots v_{p-2}v_{p-1}$. Set $\Gamma_1^*$ = $v_0$ and $\Gamma_2^*$ = $v_1v_2\cdots v_{p-2}v_{p-1}$. Let $\Lambda$ be the graph of groups obtained from $\Gamma^*$ by contracting $\Gamma_i^*$ to a point for all $i$ $\in$ \{1,2\}. Let $\Lambda'$ be the graph of groups obtained from $\Lambda$ by composing with $\pi$. As in Step 3, we may replace $\Gamma$ by $\Lambda'$, so as to assume that $\Gamma$ has two vertices and two edges joining these two vertices.
\\

\underline{End of the proof}: Present $G_v$ = $\langle g \mid g^{p^n} = 1 \rangle$, $G_w$ = $\langle h \mid h^{p^n} = 1 \rangle$, and $G$ as described above. There are two cases:
\\

\underline{Case 1}: $\Gamma$ has one edge.
\\

In this case, $G$ is an amalgamated product of two finite abelian $p$-groups. Since $G$ is residually $p$-finite, $G$ is conjugacy $p$-separable by \cite{I}, Theorem 2. Thus, there exists a homomorphism $\varphi$ from $G$ onto a finite $p$-group $P$ such that $\varphi(g)$ $\nsim$ $\varphi(h)$.
\\

\underline{Case 2}: $\Gamma$ has two edges.
\\

We have:
\begin{center}
$G$ = $\langle$ $g$, $h$, $t$ $\mid$ $g^{p^n}$ = 1, $h^{p^n}$ = 1, $g^{p^r}$ = $h^{p^r}$, $tg^{p^s}t^{-1}$ = $h^{p^s}$ $\rangle$,
\end{center}
where $r$ $\in$ \{1,...,$n-1$\}, $s$ $\in$ \{1,...,$n-1$\}. Let:
\begin{center}
$$
A = C_{p^n}\times\underbrace{C_{p^s}\times\cdots\times
C_{p^s}}_{\mbox{$p^r-1$}}\times C_{p^r}
$$
\end{center}
Set $m$ = $p^r+1$. Present each factor of this product in the natural way, using generators $x_1$,...,$x_m$ respectively. Let $\alpha$ be the automorphism of $A$ defined by:
\begin{center}
$\alpha(x_1)$ = $x_1x_2x_m$\\
$\alpha(x_i)$ = $x_{i+1}$, $\forall$ $i$ $\in$ \{2,...,$m-2$\}\\
$\alpha(x_{m-1})$ = $(x_2\cdots x_{m-1})^{-1}$\\
$\alpha(x_m)$ = $x_m$
\end{center}
It is easily seen that $\alpha$ has order $m-1$ = $p^r$. We have:
\begin{center}
$\alpha^0(x_1)$ = $x_1$,\\
$\alpha^1(x_1)$ = $x_1x_2x_m$,\\
$\alpha^2(x_1)$ = $x_1x_2x_3x_m^2$,\\
...\\
$\alpha^{m-2}(x_1)$ = $x_1x_2x_3\cdots x_{m-1}x_m^{m-2}$.
\end{center}
Let $B$ = $A\rtimes\langle\alpha\rangle$ be the semidirect product of $A$ by $\langle\alpha\rangle$. Note that $B$ is a finite $p$-group. Let $\varphi$ : $G$ $\rightarrow$ $B$ be the homomorphism defined by:
\begin{center}
$\varphi(g)$ = $x_1$,\\
$\varphi(h)$ = $x_1x_m$,\\
$\varphi(t)$ = $\alpha$.
\end{center}
Observe that the conjugacy class of $\varphi(g)$ in $B$ is $\varphi(g)^B$ = \{$\alpha^k(x_1)$ $\mid$ $k$ $\in$ \{0,...,$m-2$\}\}. Thus, $\varphi(g)$ and $\varphi(h)$ are not conjugate in $B$.\hfill$\square$

\begin{corollary}\label{4.3}
Let $P$ be a finite $p$-group. Let $A$ be a subgroup of $P$. Let $Q$ be the HNN extension of $P$ relative to $A$:
\begin{center}
$Q$ = $\langle$ $P$, $t$ $\mid$ $t^{-1}at$ = $a$, $\forall$ $a$ $\in$ $A$ $\rangle$.
\end{center}
Then $Q$ is hereditarily conjugacy $p$-separable.
\end{corollary}

\textit{Proof}: Let $R$ be an arbitrary subgroup of $Q$. Let $f$ : $Q$ $\rightarrow$ $P$ be the natural homomorphism. We have $ker(f)\cap P$ = \{1\}. Therefore $ker(f)$ is free by \cite{KS}, Theorem 6. That is, $Q$ is an extension of a free group by a finite $p$-group. Thus $R$ is itself an extension of a free group by a finite $p$-group. Therefore $R$ is conjugacy $p$-separable by Theorem \ref{4.2}.\hfill$\square$

\begin{remark}
It is known that a fundamental group of a graph of groups, whose vertex groups are all finite $p$-groups is residually $p$-finite if and only if it is an extension of a free group by a finite $p$-group (see, for example, \cite{AF1}, Lemma 3.3). Thus, as an immediate consequence of Theorem \ref{4.2}, we have that a fundamental group of a graph of groups whose vertex groups are all finite $p$-groups is conjugacy $p$-separable if and only if it is residually $p$-finite.
\end{remark}

\section{Retractions}

In this section, we shall prove several results on retractions that will allow us to control the growth of the intersection of Lemma \ref{6.5} in finite $p$-group quotients of $G_\Gamma$.

\begin{definition}
Let $G$ be a group, and $H$ be a subgroup of $G$. We say that $H$ is a retract of $G$ if there exists a homomorphism $\rho_H$ : $G$ $\rightarrow$ $G$ such that $\rho_H(G)$ = $H$ and $\rho_H(h)$ = $h$ for all $h$ $\in$ $H$. The homomorphism $\rho_H$ is called a retraction of $G$ onto $H$.
\end{definition}

\begin{remark} 
If $G$ is a right-angled Artin group, and $H$ = $\langle W \rangle$ is a special subgroup of $G$, then $H$ is a retract of $G$. A retraction of $G$ onto $H$ is given by:
\begin{center}
$\rho_H(v) = \left\{
    \begin{array}{ll}
        v & \mbox{if } v \in W \\
        1 & \mbox{if } v \in V \setminus W
    \end{array}
\right.$
\end{center}
\end{remark}

\begin{lemma}\label{5.3}
Let $G$ be a group, and $H$ be a subgroup of $G$. Suppose that $H$ is a retract of $G$. Let $\rho_H$ be a retraction of $G$ onto $H$. Let $N$ be a normal subgroup of $G$ such that $\rho_H(N)$ $\subset$ $N$. Then $\rho_H$ induces a retraction $\rho_{\overline{H}}$ : $G/N$ $\rightarrow$ $G/N$ of $G/N$ onto the canonical image $\overline{H}$ of $H$ in $G/N$, defined by: $\rho_{\overline{H}}(gN)$ = $\rho_H(g)N$ for all $gN$ $\in$ $G/N$.
\end{lemma}

\textit{Proof} : Proved in \cite{M} (see Lemma 4.1).\hfill$\square$
\\

\begin{remark}\label{5.4} 
Let $G$ be a group, and let $H$, $H'$ be two subgroups of $G$. Suppose that $H$ and $H'$ are retracts of $G$ and that the corresponding retractions, $\rho_H$ and $\rho_{H'}$, commute. Then $\rho_H(H')$ = $\rho_{H'}(H)$ = $H\cap$$H'$. Moreover $H\cap$$H'$ is a retract of $G$. A retraction of $G$ onto $H\cap$$H'$ is given by $\rho_{H\cap H'}$ = $\rho_H\circ\rho_{H'}$ = $\rho_{H'}\circ\rho_H$.
\end{remark}

\begin{proposition}\label{5.5}
Let $G$ be a group and $H_1$,...,$H_n$ be $n$ subgroups of $G$. Suppose that $H_1$,...,$H_n$ are retracts of $G$ and that the corresponding retractions pairwise commute. Then, for every normal subgroup $K$ of $p$-power index in $G$, there exists a normal subgroup $N$ of $p$-power index in $G$ such that $N$ $<$ $K$ and $\rho_{H_i}(N)$ $\subset$ $N$ for all $i$ $\in$ \{1,...,$n$\}. Consequently, for every $i$ $\in$ \{1,...,$n$\}, the retraction $\rho_{H_i}$ induces a retraction $\rho_{\overline{H_i}}$ of $G/N$ onto the canonical image $\overline{H_i}$ of $H_i$ in $G/N$.
\end{proposition}

\textit{Proof}: Proved in \cite{M} (see Proposition 4.3 and Remark 4.4).\hfill$\square$

\begin{lemma}\label{5.6}
Let $G$ be a group, and let $H$, $H'$ be two subgroups of $G$. Suppose that $H$ and $H'$ are retracts of $G$ and that the corresponding retractions, $\rho_H$ and $\rho_{H'}$, commute. Let $N$ be a normal subgroup of $G$ such that $\rho_H(N)$ $\subset$ $N$ and $\rho_{H'}(N)$ $\subset$ $N$. Then $\varphi(H\cap$$H')$ = $\varphi(H)\cap\varphi(H')$, where $\varphi$ : $G$ $\rightarrow$ $G/N$ denotes the canonical projection.
\end{lemma}

\textit{Proof}: Proved in \cite{M} (see Lemma 4.5).\hfill$\square$
\\

The next statement is the analogue of Lemma 4.6 in \cite{M}:

\begin{corollary}\label{5.7}
Let $G$ be a group and $H_1$,...,$H_n$ be $n$ subgroups of $G$. Suppose that $H_1$,...,$H_n$ are retracts of $G$ and that the corresponding retractions $\rho_{H_1}$,...,$\rho_{H_n}$ pairwise commute. Then, for every normal subgroup $K$ of $p$-power index in $G$, there exists a normal subgroup $N$ of $p$-power index in $G$ such that $N$ $<$ $K$ and $\rho_{H_i}(N)$ $\subset$ $N$, for all $i$ $\in$ \{1,...,$n$\}. Moreover, if $\varphi$ : $G$ $\rightarrow$ $G/N$ denotes the canonical projection, then $\varphi(\bigcap_{i=1}^nH_i)$ = $\bigcap_{i=1}^n\varphi(H_i)$.
\end{corollary}

\textit{Proof}: By Proposition \ref{5.5}, there exists a normal subgroup $N$ of $p$-power index in $G$ such that $N$ $<$ $K$ and $\rho_{H_i}(N)$ $\subset$ $N$ for all $i$ $\in$ \{1,...,$n$\}. We denote by $\varphi$ : $G$ $\rightarrow$ $G/N$ the canonical projection. We argue by induction on $k$ $\in$ \{1,...,$n$\} to prove that $\varphi(\bigcap_{i=1}^kH_i)$ = $\bigcap_{i=1}^k\varphi(H_i)$. If $k$ = 1, then the result is trivial. Thus we can assume that $k$ $\geq$ 2 and that the result has been proved for $k-1$. We set $H'$ = $\bigcap_{i=1}^{k-1}H_i$. By Remark \ref{5.4}, $H'$ is a retract of $G$. A retraction of $G$ onto $H'$ is given by $\rho_{H'}$ = $\rho_{H_1}\circ...\circ\rho_{H_{k-1}}$. We have:
\begin{center}
$\rho_{H'}(N)$ = $\rho_{H_1}(...(\rho_{H_{k-2}}(\rho_{H_{k-1}}(N))))$ \\
$\subset$ $\rho_{H_1}(...(\rho_{H_{k-2}}(N)))$ \\
... \\
$\subset$ $\rho_{H_1}(N)$ \\
$\subset$ $N$.
\end{center}
The retractions $\rho_{H'}$ and $\rho_{H_k}$ commute, so we can apply Lemma \ref{5.6} to conclude that $\varphi(H'\cap H_k)$ = $\varphi(H')\cap\varphi(H_k)$. By the induction hypothesis, $\varphi(H')$ = $\bigcap_{i=1}^{k-1}\varphi(H_i)$. Finally $\varphi(\bigcap_{i=1}^kH_i)$ = $\bigcap_{i=1}^k\varphi(H_i)$.\hfill$\square$
\\

In the following lemmas, $G$ is a group, and $A$, $B$ are two subgroups of $G$. We assume that $A$ and $B$ are retracts of $G$ and that the corresponding retractions, $\rho_A$ and $\rho_B$, commute.

\begin{lemma}\label{5.8}
Let $x$, $y$ $\in$ $G$. We set $\alpha$ = $\rho_A(\rho_B(x)x^{-1})x\rho_B(x^{-1})$ ($\in$ $AxB$) and $\beta$ = $\rho_A(\rho_B(y)y^{-1})y\rho_B(y^{-1})$ ($\in$ $AyB$). Then the following are equivalent:
\begin{enumerate}
\item $y$ $\in$ $AxB$,
\item $\beta$ $\in$ $\alpha^{A\cap B}$.
\end{enumerate}
\end{lemma}

\textit{Proof}: Proved in \cite{M} (see Lemma 5.1).\hfill$\square$

\begin{lemma}\label{5.9}
Let $x$ $\in$ $G$. We set $\alpha$ = $\rho_A(\rho_B(x)x^{-1})x\rho_B(x^{-1})$ ($\in$ $AxB$) and $\gamma$ = $\rho_A(\rho_B(x)x^{-1})$ ($\in$ $A$). Then we have:
\begin{center}
$A\cap$$xBx^{-1}$ = $\gamma^{-1}C_{A\cap B}(\alpha)\gamma$.
\end{center}
\end{lemma}

\textit{Proof}: Proved in \cite{M} (see Lemma 5.2).\hfill$\square$
\\

The next five statements are the analogues of some statements in \cite{M} (Lemma 5.3, Corollary 5.4, Lemma 5.5, Lemma 5.6, and Lemma 5.7 respectively):

\begin{lemma}\label{5.10}
Let $x$ $\in$ $G$. We set: $\alpha$ = $\rho_A(\rho_B(x)x^{-1})x\rho_B(x^{-1})$ ($\in$ $AxB$). If $\alpha^{A\cap B}$ is finitely $p$-separable in $G$, then $AxB$ is also finitely $p$-separable in $G$.
\end{lemma}

\textit{Proof}: Let $y$ $\in$ $G$ such that $y$ $\notin$ $AxB$. We set $\beta$ = $\rho_A(\rho_B(y)y^{-1})y\rho_B(y^{-1})$. By Lemma \ref{5.8}, we have $\beta$ $\notin$ $\alpha^{A\cap B}$. Since $\alpha^{A\cap B}$ is finitely $p$-separable in $G$, there exists a normal subgroup $K$ of $p$-power index in $G$ such that, if $\psi$ : $G$ $\rightarrow$ $G/K$ denotes the canonical projection, we have: $\psi(\beta)$ $\notin$ $\psi(\alpha^{A\cap B})$ = $\psi(\alpha)^{\psi(A\cap B)}$. By Corollary \ref{5.7}, there exists a normal subgroup $N$ of $p$-power index in $G$ such that $N$ $<$ $K$, $\rho_A(N)$ $\subset$ $N$, $\rho_B(N)$ $\subset$ $N$ and, if $\varphi$ : $G$ $\rightarrow$ $G/N$ denotes the canonical projection, then: $\varphi(A\cap B)$ = $\varphi(A)\cap\varphi(B)$. Assume that $\varphi(\beta)$ $\in$ $\varphi(\alpha)^{\varphi(A\cap B)}$. Let $g$ $\in$ $A\cap B$ be such that $\varphi(\beta)$ = $\varphi(g)\varphi(\alpha)\varphi(g)^{-1}$. Then $\beta$ $\in$ $g\alpha g^{-1}N$. Since $N$ $<$ $K$, we obtain $\beta$ $\in$ $g\alpha g^{-1}K$. But this contradicts the fact that $\psi(\beta)$ $\notin$ $\psi(\alpha)^{\psi(A\cap B)}$. Therefore we have: $\varphi$($\beta$) $\notin$ $\varphi(\alpha)^{\varphi(A\cap B)}$ i.e. $\varphi(\beta)$ $\notin$ $\varphi(\alpha)^{\varphi(A)\cap\varphi(B)}$. We set $\overline{A}$ = $\varphi(A)$ and $\overline{B}$ = $\varphi(B)$. By Lemma \ref{5.3}, $\rho_A$ induces a retraction $\rho_{\overline{A}}$ of $G/N$ onto $\overline{A}$ and $\rho_B$ induces a retraction $\rho_{\overline{B}}$ of $G/N$ onto $\overline{B}$. We set: $\overline{x}$ = $\varphi(x)$ and $\overline{y}$ = $\varphi(y)$. We have: $\varphi(\alpha)$ = $\rho_{\overline{A}}(\rho_{\overline{B}}(\overline{x})\overline{x}^{-1})\overline{x}\rho_{\overline{B}}(\overline{x}^{-1})$ and $\varphi(\beta)$ = $\rho_{\overline{A}}(\rho_{\overline{B}}(\overline{y})\overline{y}^{-1})\overline{y}\rho_{\overline{B}}(\overline{y}^{-1})$. By Lemma \ref{5.8}, we have $\overline{y}$ $\notin$ $\overline{A}\overline{x}\overline{B}$ i.e. $\varphi(y)$ $\notin$ $\varphi(AxB)$.\hfill$\square$

\begin{corollary}\label{5.11}
Let $G$ be a group, and $A$, $B$ be two subgroups of $G$. Suppose that $G$ is residually $p$-finite. If $A$ and $B$ are retracts of $G$, such that the corresponding retractions commute, then $AB$ is finitely $p$-separable in $G$.
\end{corollary}

\textit{Proof}: We apply Lemma \ref{5.10} to $x$ = 1.\hfill$\square$

\begin{lemma}\label{5.12}
Let $G$ be a group, and $A$ be a subgroup of $G$. Suppose that $G$ is residually $p$-finite and that $A$ is a retract of $G$. Then if a subset $S$ of $A$ is closed in the pro-$p$ topology on $A$, it is also closed in the pro-$p$ topology on $G$.
\end{lemma}

\textit{Proof}: We denote by $\overline{S}$ the closure of $S$ in $G$ -- equipped with the pro-$p$ topology. We shall show that $\overline{S}$ $\subset$ $S$. By Corollary \ref{5.11}, $A$ is closed in $G$. Therefore $\overline{S}$ $\subset$ $A$. Let $a$ $\in$ $G\setminus$$S$. We can assume that $a$ $\in$ $A$. There exists a homomorphism $\psi$ from $A$ onto a finite $p$-group $P$ such that $\psi(a)$ $\notin$ $\psi(S)$. We set: $\varphi$ = $\psi\circ\rho_A$. We have: $\varphi(a)$ = $\psi(a)$ $\notin$ $\psi(S)$ = $\varphi(S)$. Then $a$ $\notin$ $\overline{S}$.\hfill$\square$

\begin{lemma}\label{5.13}
Let $x$ $\in$ $G$. We set $\alpha$ = $\rho_A(\rho_B(x)x^{-1})x\rho_B(x^{-1})$. Suppose that the pair ($A\cap$$B$,$\alpha$) satisfies the $p$-centralizer condition in $G$. Then, for every normal subgroup $K$ of $p$-power index in $G$, there exists a normal subgroup $N$ of $p$-power index in $G$ such that $N$ $<$ $K$, $\rho_A(N)$ $\subset$ $N$, $\rho_B(N)$ $\subset$ $N$ and, if $\varphi$ : $G$ $\rightarrow$ $G/N$ denotes the canonical projection, $\varphi(A)\cap\varphi(xBx^{-1})$ $\subset$ $\varphi(A\cap xBx^{-1})\varphi(K)$.
\end{lemma}

\textit{Proof}: Let $K$ be a normal subgroup of $p$-power index in $G$. We set $\gamma$ = $\rho_A(\rho_B(x)x^{-1})$ $\in$ $A$. By Lemma \ref{5.9}, we have: $A\cap$$xBx^{-1}$ = $\gamma^{-1}C_{A\cap B}(\alpha)\gamma$. Since the pair ($A\cap B$,$\alpha$) satisfies $pCC_G$, there exists a normal subgroup $L$ of $p$-power index in $G$ such that $L$ $<$ $K$ and, if $\psi$ : $G$ $\rightarrow$ $G/L$ denotes the canonical projection, $C_{\psi(A\cap B)}(\psi(\alpha))$ $\subset$ $\psi(C_{A\cap B}(\alpha)K)$. This is equivalent to $\psi^{-1}(C_{\psi(A\cap B)}(\psi(\alpha)))$ $\subset$ $C_{A\cap B}(\alpha)K$. Indeed let $g$ $\in$ $\psi^{-1}(C_{\psi(A\cap B)}(\psi(\alpha)))$. We have $\psi(g)$ $\in$ $C_{\psi(A\cap B)}(\psi(\alpha))$ $\subset$ $\psi(C_{A\cap B}(\alpha)K)$. Then $g$ $\in$ $C_{A\cap B}(\alpha)KL$ $\subset$ $C_{A\cap B}(\alpha)K$ (because $L$ $<$ $K$). By Corollary \ref{5.7}, there exists a normal subgroup $N$ of $p$-power index in $G$ such that $N$ $<$ $L$, $\rho_A(N)$ $\subset$ $N$, $\rho_B(N)$ $\subset$ $N$ and, if $\varphi$ : $G$ $\rightarrow$ $G/N$ denotes the canonical projection, $\varphi(A\cap B)$ = $\varphi(A)\cap\varphi(B)$. We set $\overline{A}$ = $\varphi(A)$, $\overline{B}$ = $\varphi(B)$. By Lemma \ref{5.3}, $\rho_A$ induces a retraction $\rho_{\overline{A}}$ of $G/N$ onto $\overline{A}$, and $\rho_B$ induces a retraction $\rho_{\overline{B}}$ of $G/N$ onto $\overline{B}$. Obviously $\rho_{\overline{A}}$ and $\rho_{\overline{B}}$ commute. We set $\overline{x}$ = $\varphi(x)$, $\overline{\alpha}$ = $\rho_{\overline{A}}(\rho_{\overline{B}}(\overline{x})\overline{x}^{-1})\overline{x}\rho_{\overline{B}}(\overline{x}^{-1})$ ($\in$ $G/N$) and $\overline{\gamma}$ = $\rho_{\overline{A}}(\rho_{\overline{B}}(\overline{x})\overline{x}^{-1})$ ($\in$ $\overline{A}$). Observe that $\overline{\alpha}$ = $\varphi(\alpha)$ and $\overline{\gamma}$ = $\varphi(\gamma)$. Then, by Lemma \ref{5.9}, we have: $\overline{A}\cap$$\overline{x}\overline{B}\overline{x}^{-1}$ = $\overline{\gamma}^{-1}C_{\overline{A}\cap\overline{B}}(\overline{\alpha})\overline{\gamma}$. Now, $\overline{A}\cap\overline{B}$ = $\varphi(A\cap B)$. Thus:
\begin{center}
$\varphi^{-1}(\overline{A}\cap\overline{x}\overline{B}\overline{x}^{-1})$ = $\varphi^{-1}(\overline{\gamma}^{-1}C_{\varphi(A\cap B)}(\overline{\alpha})\overline{\gamma})$ = $\gamma^{-1}\varphi^{-1}(C_{\varphi(A\cap B)}(\overline{\alpha}))\gamma$.
\end{center}
We claim that:
\begin{center}
$\varphi^{-1}(C_{\varphi(A\cap B)}(\varphi(\alpha)))$ $\subset$ $\psi^{-1}(C_{\psi(A\cap B)}(\psi(\alpha)))$.
\end{center}
Indeed let $g$ $\in$ $\varphi^{-1}(C_{\varphi(A\cap B)}(\varphi(\alpha)))$. We have $\varphi(g)$ $\in$ $\varphi(A\cap B)$ i.e. $g$ $\in$ $(A\cap B)N$, which implies $g$ $\in$ $(A\cap B)L$ i.e. $\psi(g)$ $\in$ $\psi(A\cap B)$; and $\varphi(g)\varphi(\alpha)$ = $\varphi(\alpha)\varphi(g)$ i.e. $g\alpha g^{-1}\alpha^{-1}$ $\in$ $N$, which implies $g\alpha g^{-1}\alpha^{-1}$ $\in$ $L$ i.e. $\psi(g)\psi(\alpha)$ = $\psi(\alpha)\psi(g)$. We deduce that:
\begin{center}
$\varphi^{-1}(C_{\varphi(A\cap B)}(\varphi(\alpha)))$ $\subset$ $C_{A\cap B}(\alpha)K$,
\end{center}
and hence:
\begin{center}
$\varphi^{-1}(\overline{A}\cap\overline{x}\overline{B}\overline{x}^{-1})$ $\subset$ $\gamma^{-1}C_{A\cap B}(\alpha)\gamma K$ = $(A\cap$$xBx^{-1})K$.
\end{center}
We conclude that:
\begin{center}
$\varphi(A)\cap\varphi(xBx^{-1})$ $\subset$ $\varphi(A\cap xBx^{-1})\varphi(K)$.
\end{center}\hfill$\square$

\begin{lemma}\label{5.14}
Let $x$, $y$ $\in$ $G$. We set $C$ = $xBx^{-1}$ ($<$ $G$) and $\alpha$ = $\rho_A(\rho_B(x)\\x^{-1})x\rho_B(x^{-1})$. If $\alpha^{A\cap B}$ and $y^{A\cap C}$ are finitely $p$-separable in $G$, and if the pair ($A\cap B$, $\alpha$) satisfies $pCC_G$, then $C_A(y)C$ is finitely $p$-separable in $G$.
\end{lemma}

\textit{Proof}: Let $z$ $\in$ $G$ such that $z$ $\notin$ $C_A(y)C$. Suppose first that $z$ $\notin$ $AC$. Since $\alpha^{A\cap B}$ is finitely $p$-separable in $G$, $AxB$ is finitely $p$-separable in $G$ by Lemma \ref{5.10}. Therefore $AC$ = $AxBx^{-1}$ is also finitely $p$-separable in $G$. Consequently there exists a normal subgroup $N$ of $p$-power index in $G$ such that $z$ $\notin$ $ACN$. We obviously have $z$ $\notin$ $C_A(y)CN$. Thus we can assume that $z$ $\in$ $AC$. Let $a$ $\in$ $A$, $c$ $\in$ $C$ be such that $z$ = $ac$. Since $z$ $\notin$ $C_A(y)C$, $a^{-1}ya$ $\notin$ $y^{A\cap C}$. Indeed, if there is $g$ $\in$ $A\cap C$ such that $a^{-1}ya$ = $gyg^{-1}$, then $(ag)^{-1}y(ag)$ = $y$ i.e. $ag$ $\in$ $C_A(y)$. We obtain $a$ $\in$ $C_A(y)C$, and then $z$ $\in$ $C_A(y)C$ -- which is a contradiction. Now $y^{A\cap C}$ is finitely $p$-separable in $G$. Then there exists a normal subgroup $K$ of $p$-power index in $G$ such that $a^{-1}ya$ $\notin$ $y^{A\cap C}K$. By Lemma \ref{5.13}, there exists a normal subgroup $N$ of $p$-power index in $G$ such that $N$ $<$ $K$ and, if $\varphi$ : $G$ $\rightarrow$ $G/N$ denotes the canonical projection, $\varphi(A)\cap\varphi(C)$ $\subset$ $\varphi(A\cap C)\varphi(K)$. For a subset $S$ of $G$, we set $\overline{S}$ = $\varphi(S)$. For an element $g$ of $G$, we set $\overline{g}$ = $\varphi(g)$. We have: $\overline{y}^{\overline{A}\cap \overline{C}}$ $\subset$ $\overline{y}^{\overline{A\cap C}.\overline{K}}$. Note that $\overline{K}$ $\lhd$ $G/N$. Then $\overline{y}^{\overline{A}\cap \overline{C}}$ $\subset$ $\overline{y}^{\overline{A\cap C}}\overline{K}$. Observe that $\overline{a}^{-1}\overline{y}$\hspace{1mm}$\overline{a}$ $\notin$ $\overline{y}^{\overline{A\cap C}}\overline{K}$ -- otherwise we would have $a^{-1}ya$ $\in$ $y^{A\cap C}KN$, and then $a^{-1}ya$ $\in$ $y^{A\cap C}K$ (because $N$ $<$ $K$). We deduce that $\overline{a}^{-1}\overline{y}$\hspace{1mm}$\overline{a}$ $\notin$ $\overline{y}^{\overline{A}\cap \overline{C}}$. Now it suffices to show that $\varphi(z)$ $\notin$ $\varphi(C_A(y)C)$. Suppose the contrary. Let $a'$ $\in$ $C_A(y)$, $c'$ $\in$ $C$ be such that $\varphi(z)$ = $\varphi(a'c')$. Then $\varphi(ac)$ = $\varphi(a'c')$. Thus $\varphi(a'^{-1}a)$ = $\varphi(c'c^{-1})$. We set $\overline{g}$ = $\varphi(a'^{-1}a)$ = $\varphi(c'c^{-1})$ ($\in$ $\overline{A}\cap\overline{C}$). We have: $\varphi(z)$ = $\varphi(a')\overline{g}\varphi(c)$ and $\overline{a}$ = $\varphi(z)\varphi(c)^{-1}$ = $\varphi(a')\overline{g}$. Then $\overline{a}^{-1}\overline{y}$\hspace{1mm}$\overline{a}$ = $\overline{g}^{-1}\varphi(a'^{-1}ya')\overline{g}$ = $\overline{g}^{-1}\varphi(y)\overline{g}$ = $\overline{g}^{-1}\overline{y}\hspace{1mm}\overline{g}$ $\in$ $\overline{y}^{\overline{A}\cap\overline{C}}$ -- a contradiction. We have shown that $C_A(y)C$ is finitely $p$-separable in $G$.\hfill$\square$

\section{Proof of the main theorem}

\hspace{5mm}We turn now to the proof that right-angled Artin groups are hereditarily conjugacy $p$-separable. We need the following theorem, which is due to Duchamp and Krob (see \cite{DK2}, Theorem 2.3).

\begin{theorem}\label{6.1}
Right-angled Artin groups are residually $p$-finite.
\end{theorem}

(Note that the exact statement of \cite{DK2}, Theorem 2.3, is that right-angled Artin groups are residually torsion-free nilpotent; Theorem \ref{6.1} then follows from \cite{G}, Theorem 2.1.) This theorem can also be proved using HNN extensions (see \cite{Lo}, Theorem 2.11).
\\

Basically, Proposition \ref{6.2} establishes the main result. Proposition \hyperlink{6.2.1}{6.2.1} and Proposition \hyperlink{6.2.2}{6.2.2} will be proved simultaneously by induction on the rank of $G$.

\begin{proposition}\label{6.2}
Let $G$ be a right-angled Artin group.
\begin{description}
\item\hypertarget{6.2.1}{6.2.1} Every special subgroup $S$ of $G$ satisfies the $p$-centralizer condition in $G$ ($pCC_G$).
\item\hypertarget{6.2.2}{6.2.2} For all $g$ $\in$ $G$ and for every special subgroup $S$ of $G$, $g^S$ is finitely $p$-separable in $G$.
\end{description}
\end{proposition}

From now on, we assume that $G$ is a right-angled Artin group of rank $r$ ($r$ $\geq$ 1), and that $H$ = $\langle W \rangle$ is a special subgroup of $G$ of rank $r-1$. Thus, $G$ can be written as an HNN extension of $H$ relative to the special subgroup $K$ = $\langle link(t) \rangle$ of $H$:
\begin{center}
$G$ = $\langle$ $H$, $t$ $\mid$ $t^{-1}kt$ = $k$, $\forall$ $k$ $\in$ $K$ $\rangle$.
\end{center}
Recall that $H$ is a retract of $G$. A retraction of $G$ onto $H$ is given by:
$$
\rho_H(v) = \left\{
    \begin{array}{ll}
        v & \mbox{if } v \in W \\
        1 & \mbox{if } v \in V \setminus W
    \end{array}
\right.
$$
We also assume that:
\begin{itemize}
\item every special subgroup $S$ of $H$ satisfies the $p$-centralizer condition in $H$ ($pCC_H$),
\item for all $h$ $\in$ $H$ and for every special subgroup $S$ of $H$, $h^S$ is finitely $p$-separable in $H$.
\end{itemize}

The next results (Lemma \ref{6.3} to Lemma \ref{6.14}) are preliminaries to the proof of Proposition \ref{6.2}.
\\

In general, if $A$ and $B$ are subgroups of a group $G$, the image of the intersection of $A$ and $B$ under a homomorphism $\varphi$ : $G$ $\rightarrow$ $H$ does not coincide with the intersection of the images of $A$ and $B$ in $H$. However, the $p$-centralizer condition and the above results on retractions will allow us to obtain the following lemma, which will be used to apply Minasyan's criterion for conjugacy in HNN extensions (see Lemma \ref{6.5}).

\begin{lemma}\label{6.3}
Let be given $A_0$, a conjugate of a special subgroup of $H$, $A_1$,...,$A_n$, $n$ special subgroups of $H$, and $\alpha$,$x_0$,...,$x_n$,$y_1$,...,$y_n$, 2($n$+1) elements of $H$. Then, for every normal subgroup $L$ of $p$-power index in $H$, there exists a normal subgroup $N$ of $p$-power index in $H$ such that $N$ $<$ $L$ and, if $\varphi$ : $H$ $\rightarrow$ $H/N$ denotes the canonical projection, then:
\begin{center}
$\overline{\alpha}C_{\overline{A_0}}(\overline{x_0})\cap\bigcap_{i=1}^n\overline{x_i}\overline{A_i}\overline{y_i}$ $\subset$ $\varphi((\alpha C_{A_0}(x_0)\cap\bigcap_{i=1}^nx_iA_iy_i)L)$,
\end{center}
where $\overline{A_i}$ = $\varphi(A_i)$ ($i$ $\in$ \{0,...,$n$\}), $\overline{\alpha}$ = $\varphi(\alpha)$, $\overline{x_j}$ = $\varphi(x_j)$ ($j$ $\in$ \{0,...,$n$\}), $\overline{y_k}$ = $\varphi(y_k)$ ($k$ $\in$ \{1,...,$n$\}).
\end{lemma}

\textit{Proof}: Let $L$ be a subgroup of $p$-power index in $H$. We argue by induction on $n$. Strictly speaking, the basis of our induction is n = 0 but we will need the case n = 1. By the assumptions, there exist a special subgroup $A$ of $H$, and an element $\beta$ of $H$ such that $A_0$ = $\beta A \beta^{-1}$. 
\\

\underline{n = 0}: We set $x$ = $\beta^{-1}x_0\beta$. The pair ($A$,$x$) satisfies $pCC_H$ by the assumptions. There exists a normal subgroup $N$ of $p$-power index in $H$ such that $N$ $<$ $L$ and, if $\varphi$ : $H$ $\rightarrow$ $H/N$ denotes the canonical projection, then $C_{\varphi(A)}(\varphi(x))$ $\subset$ $\varphi(C_A(x)L)$. But $C_{A_0}(x_0)$ = $\beta C_A(x) \beta^{-1}$. We deduce that: $\varphi(\alpha)C_{\varphi(A_0)}(\varphi(x_0))$ $\subset$ $\varphi((\alpha C_{A_0}(x_0))L)$. 
\\

\underline{n = 1}: There are two cases:

\vspace{1mm}

\underline{Case 1}: $\alpha C_{A_0}(x_0) \cap x_1A_1y_1$ = $\emptyset$. This is equivalent to say that: $x_1$ $\notin$ $\alpha C_{A_0}(x_0)y_1^{-1}A_1$. We set $B$ = $(y_1\beta)^{-1}A_1y_1\beta$, so that we have: $x_1$ $\notin$ $\alpha\beta(C_A(x)B)\beta^{-1}y_1^{-1}$. Now the intersection of conjugates of two special subgroups of $H$ is a conjugate of a special subgroup of $H$ (see \cite{M}, Lemma 6.5). Then $A\cap A_1$ is a conjugate of a special subgroup $C$ of $H$. There exists $\gamma$ $\in$ $H$ such that $A \cap A_1$ = $\gamma C \gamma^{-1}$. Therefore if $h$ $\in$ $H$, $h^{A \cap A_1}$ = $\gamma(\gamma^{-1}h\gamma)^C\gamma^{-1}$. Now $(\gamma^{-1}h\gamma)^C$ is finitely $p$-separable in $H$ by the assumptions. We deduce that $h^{A \cap A_1}$ is finitely $p$-separable in $H$. With the same argument, $x^{A \cap B}$ is finitely $p$-separable in $H$. Now the pair ($A\cap A_1$,$h$) satisfies $pCC_H$ by the assumptions. We deduce that $C_A(x)B$ is finitely $p$-separable in $H$ by Lemma \ref{5.14}. This implies that $\alpha C_{A_0}(x_0)y_1^{-1}A_1$ is finitely $p$-separable in $H$. There exists a normal subgroup $M$ of $p$-power index in $H$ such that $x_1$ $\notin$ $\alpha C_{A_0}(x_0)y_1^{-1}A_1M$. Up to replacing $M$ by $M\cap L$, we can assume that $M$ $<$ $L$. Now the pair ($A_0$,$x_0$) satisfies $pCC_H$ by the assumptions. There exists a normal subgroup $N$ of $p$-power index in $H$ such that $N$ $<$ $M$ and, if $\varphi$ : $H$ $\rightarrow$ $H/N$ denotes the canonical projection, then $C_{\varphi(A_0)}(\varphi(x_0))$ $\subset$ $\varphi(C_{A_0}(x_0)M)$, or, equivalently, $\varphi^{-1}(C_{\varphi(A_0)}(\varphi(x_0)))$ $\subset$ $C_{A_0}(x_0)M$. Then $\varphi^{-1}(\overline{\alpha}C_{\overline{A_0}}(\overline{x_0})\overline{y_1}^{-1}\overline{A_1})$ $\subset$ $\alpha\varphi^{-1}(C_{\overline{A_0}}(\overline{x_0}))y_1^{-1}A_1$ $\subset$ $\alpha C_{A_0}(x_0)y_1^{-1}A_1M$ (with the same notations as in the statement of the lemma). Therefore: $x_1$ $\notin$ $\varphi^{-1}(\overline{\alpha}C_{\overline{A_0}}(\overline{x_0})\overline{y_1}^{-1}\overline{A_1})$. Finally: $\overline{\alpha}C_{\overline{A_0}}(\overline{x_0})\cap \overline{x_1}\overline{A_1}\overline{y_1}$ = $\emptyset$.

\vspace{1mm}

\underline{Case 2}: $\alpha C_{A_0}(x_0) \cap x_1A_1y_1$ $\neq$ $\emptyset$.
\begin{remark}\label{6.4}
If $G$ is a group, and $H$, $K$ are two subgroups of $G$ such that $aH\cap bKc$ $\neq$ $\emptyset$ -- where $a$, $b$, $c$ $\in$ $G$ --, then for all $g$ $\in$ $aH\cap bKc$, we have $aH\cap bKc$ = $g(H\cap c^{-1}Kc)$. 
\end{remark}
Choose $g$ $\in$ $\alpha C_{A_0}(x_0) \cap x_1A_1y_1$. By Remark \ref{6.4}, we have: $\alpha C_{A_0}(x_0) \cap x_1A_1y_1$ = $g(C_{A_0}(x_0)\cap y_1^{-1}A_1y_1)$. We set $D$ = $A_0\cap y_1^{-1}A_1y_1$. Then $\alpha C_{A_0}(x_0) \cap x_1A_1y_1$ = $gC_D(x_0)$. Now, $D$ is a conjugate of a special subgroup $E$ of $H$ by \cite{M}, Lemma 6.5. There exists $\delta$ $\in$ $H$ such that $D$ = $\delta E\delta^{-1}$. As above, the pair ($D$,$x_0$) satisfies $pCC_H$. There exists a normal subgroup $M$ of $p$-power index in $H$ such that $M$ $<$ $L$ and, if $\psi$ : $H$ $\rightarrow$ $H/M$ denotes the canonical projection, we have: $C_{\psi(D)}(\psi(x_0))$ $\subset$ $\psi(C_D(x_0)L)$. Now by Lemma \ref{5.13}, there exists a normal subgroup $N$ of $p$-power index in $H$ such that $N$ $<$ $M$ and, if $\varphi$ : $H$ $\rightarrow$ $H/N$ denotes the canonical projection, then $\varphi(A)\cap\varphi((y_1\beta)^{-1}A_1y_1\beta)$ $\subset$ $\varphi(A\cap (y_1\beta)^{-1}A_1y_1\beta)\varphi(M)$. Therefore:
\begin{center}
$\overline{A_0}\cap \overline{y_1}^{-1}\overline{A_1}\overline{y_1}$ = $\varphi(\beta A\beta^{-1})\cap\varphi(y_1^{-1}A_1y_1)$ = $\varphi(\beta)(\varphi(A)\cap\varphi((y_1\beta)^{-1}A_1y_1\beta))\varphi(\beta^{-1})$ $\subset$ $\varphi(\beta)(\varphi(A\cap(y_1\beta)^{-1}A_1y_1\beta)\varphi(M))\varphi(\beta^{-1})$ = $\varphi(A_0\cap y_1^{-1}A_1y_1)\varphi(M)$ = $\varphi(D)\varphi(M)$ ($\ast$),
\end{center}
(with the same notations as in the statement of the lemma). We set $\overline{g}$ = $\varphi(g)$. Note that $\overline{g}$ $\in$ $\overline{\alpha} C_{\overline{A_0}}(\overline{x_0})\cap \overline{x_1}\overline{A_1}\overline{y_1}$. Therefore $\overline{\alpha} C_{\overline{A_0}}(\overline{x_0})\cap \overline{x_1}\overline{A_1}\overline{y_1}$ = $\overline{g}(C_{\overline{A_0}}(\overline{x_0})\cap \overline{y_1}^{-1}\overline{A_1}\overline{y_1})$. Considering ($\ast$), we obtain:
\begin{center}
$\overline{\alpha} C_{\overline{A_0}}(\overline{x_0})\cap \overline{x_1}\overline{A_1}\overline{y_1}$ = $\overline{g}C_{\overline{A_0}\cap\overline{y_1}^{-1}\overline{A_1}\overline{y_1}}(\overline{x_0})$ $\subset$ $\overline{g}C_{\varphi(D)\varphi(M)}(\overline{x_0})$.
\end{center}
Recall that $N$ $<$ $M$. Then $\psi$ : $H$ $\rightarrow$ $H/M$ induces a homomorphism $\widetilde{\psi}$ : $H/N$ $\rightarrow$ $H/M$ such that $\psi$ = $\widetilde{\psi}\circ\varphi$. Note that $\widetilde{\psi}(\varphi(D)\varphi(M))$ = $\psi(D)$. Let $z$ $\in$ $C_{\varphi(D)\varphi(M)}(\overline{x_0})$. Then:
\begin{center}
$\widetilde{\psi}(z)$ $\in$ $C_{\psi(D)}(\psi(x_0))$ $\subset$ $\psi(C_D(x_0)L)$ = $\widetilde{\psi}(\varphi(C_D(x_0)L))$.
\end{center}
Therefore $z$ $\in$ $\varphi(C_D(x_0)L)ker(\widetilde{\psi})$ = $\varphi(C_D(x_0)L)$ (because $ker(\widetilde{\psi})$ = $\varphi(M)$ $<$ $\varphi(L)$). We deduce that $C_{\varphi(D)\varphi(M)}(\overline{x_0})$ $\subset$ $\varphi(C_D(x_0)L)$. We conclude that
\begin{center}
$\overline{\alpha} C_{\overline{A_0}}(\overline{x_0})\cap \overline{x_1}\overline{A_1}\overline{y_1}$ $\subset$ $\overline{g}\varphi(C_D(x_0)L)$ = $\varphi(gC_D(x_0)L)$ = $\varphi((\alpha C_{A_0}(x_0)\cap x_1A_1y_1)L)$.
\end{center}

\vspace{2mm}

\underline{Inductive step}: Suppose that $n$ $\geq$ 1 and that the result has been proved for $n-1$. Note that if $\alpha C_{A_0}(x_0)\cap\bigcap_{i=1}^{n-1}x_iA_iy_i$ = $\emptyset$, then by the induction hypothesis, there exists a normal subgroup $N$ of $p$-power index in $H$ such that, if $\varphi$ : $H$ $\rightarrow$ $H/N$ denotes the canonical projection, then
\begin{center}
$\overline{\alpha}C_{\overline{A_0}}(\overline{x_0})\cap\bigcap_{i=1}^{n-1}\overline{x_i}\overline{A_i}\overline{y_i}$ $\subset$ $\varphi((\alpha C_{A_0}(x_0)\cap\bigcap_{i=1}^{n-1}x_iA_iy_i)L)$ = $\emptyset$.
\end{center}
Obviously:
\begin{center}
$\overline{\alpha}C_{\overline{A_0}}(\overline{x_0})\cap\bigcap_{i=1}^n\overline{x_i}\overline{A_i}\overline{y_i}$ = $\emptyset$ $\subset$ $\varphi((\alpha C_{A_0}(x_0)\cap\bigcap_{i=1}^nx_iA_iy_i)L)$
\end{center}
Thus we can assume that $\alpha C_{A_0}(x_0)\cap\bigcap_{i=1}^{n-1}x_iA_iy_i$ $\neq$ $\emptyset$. Therefore $\alpha C_{A_0}(x_0)\cap\bigcap_{i=1}^{n-1}x_iA_iy_i$ = $g$($C_{A_0}(x_0)\cap\bigcap_{i=1}^{n-1}y_i^{-1}A_iy_i$) for some $g$ $\in$ $H$. We set $F$ = $A_0\cap\bigcap_{i=1}^{n-1}y_i^{-1}A_iy_i$. Again, $F$ is a conjugate of a special subgroup of $H$ by \cite{M}, Lemma 6.5. We have: $\alpha C_{A_0}(x_0)\cap\bigcap_{i=1}^{n-1}x_iA_iy_i$ = $gC_F(x_0)$. Now, by the case n = 1, there exists a normal subgroup $M$ of $p$-power index in $H$ such that $M$ $<$ $L$ and, if $\psi$ : $H$ $\rightarrow$ $H/M$ denotes the canonical projection, then:
\begin{center}
$\psi(g)C_{\psi(F)}(\psi(x_0))\cap\psi(x_nA_ny_n)$ $\subset$ $\psi((gC_F(x_0)\cap x_nA_ny_n)L)$,
\end{center}
or, equivalently:
\begin{center}
$\psi^{-1}(\psi(g)C_{\psi(F)}(\psi(x_0))\cap\psi(x_nA_ny_n))$ $\subset$ $(gC_F(x_0)\cap x_nA_ny_n)L$.
\end{center}
On the other hand, by the induction hypothesis, there exists a normal subgroup $N$ of $p$-power index in $H$ such that $N$ $<$ $M$ and, if $\varphi$ : $H$ $\rightarrow$ $H/N$ denotes the canonical projection, then:
\begin{center}
$\overline{\alpha}C_{\overline{A_0}}(\overline{x_0})\cap\bigcap_{i=1}^{n-1}\overline{x_i}\overline{A_i}\overline{y_i}$ $\subset$ $\varphi((\alpha C_{A_0}(x_0)\cap\bigcap_{i=1}^{n-1}x_iA_iy_i)M)$
\end{center}
or, equivalently:
\begin{center}
$\varphi^{-1}(\overline{\alpha}C_{\overline{A_0}}(\overline{x_0})\cap\bigcap_{i=1}^{n-1}\overline{x_i}\overline{A_i}\overline{y_i})$ $\subset$ $(\alpha C_{A_0}(x_0)\cap\bigcap_{i=1}^{n-1}x_iA_iy_i)M$
\end{center}
Thus we have:
\begin{center}
$\varphi^{-1}(\overline{\alpha}C_{\overline{A_0}}(\overline{x_0})\cap\bigcap_{i=1}^n\overline{x_i}\overline{A_i}\overline{y_i})$ = $\varphi^{-1}(\overline{\alpha}C_{\overline{A_0}}(\overline{x_0})\cap\bigcap_{i=1}^{n-1}\overline{x_i}\overline{A_i}\overline{y_i})\cap \varphi^{-1}(\overline{x_n}\overline{A_n}\overline{y_n})$ $\subset$ $(\alpha C_{A_0}(x_0)\cap\bigcap_{i=1}^{n-1}x_iA_iy_i)M$$\cap$$x_nA_ny_nN$ = $gC_F(x_0)M\cap x_nA_ny_nN$.
\end{center}
Recall that $N$ $<$ $M$. Finally we have:
\begin{center}
$\varphi^{-1}(\overline{\alpha}C_{\overline{A_0}}(\overline{x_0})\cap\bigcap_{i=1}^n\overline{x_i}\overline{A_i}\overline{y_i})$ $\subset$ $gC_F(x_0)M$$\cap$$x_nA_ny_nM$ $\subset$ $\psi^{-1}(\psi(g)C_{\psi(F)}(\psi(x_0)))\cap\psi^{-1}(\psi(x_nA_ny_n))$ = $\psi^{-1}(\psi(g)C_{\psi(F)}(\psi(x_0))\cap\psi(x_nA_ny_n))$ $\subset$ $(gC_F(x_0)\cap x_nA_ny_n)L$ = $(\alpha C_{A_0}(x_0)\cap\bigcap_{i=1}^nx_iA_iy_i)L$.
\end{center}
\hfill$\square$
\\

We need the following criterion for conjugacy in HNN extensions:

\begin{lemma}\label{6.5}
Let $G$ = $\langle$ $H$, $t$ $\mid$ $t^{-1}kt$ = $k$, $\forall$ $k$ $\in$ $K$ $\rangle$ be an HNN extension. Let $S$ be a subgroup of $H$. Let $g$ = $x_0t^{a_1}x_1\cdots t^{a_n}x_n$ ($n$ $\geq$ 1) and $h$ = $y_0t^{b_1}y_1\cdots t^{b_m}y_m$ be elements of $G$ in reduced form. Then $h$ $\in$ $g^S$ if and only if all of the following conditions hold:
\begin{enumerate}
\item $m$ = $n$ and $a_i$ = $b_i$, for all $i$ $\in$ \{1,...,$n$\},
\item $y_0\cdots y_n$ $\in$ $(x_0\cdots x_n)^S$,
\item if $\alpha$ $\in$ $S$ satisfies $y_0\cdots y_n$ = $\alpha x_0\cdots x_n\alpha^{-1}$, then:
\begin{center}
$\alpha C_S(x_0\cdots x_n)\cap y_0Kx_0^{-1}\cap (y_0y_1)K(x_0x_1)^{-1}\cap\cdots \cap(y_0\cdots y_{n-1})K(x_0\cdots x_{n-1})^{-1}$ $\neq$ $\emptyset$.
\end{center}
\end{enumerate}
\end{lemma}

\textit{Proof}: Proved in \cite{M} (see Lemma 7.11).\hfill$\square$
\\

The following is the analogue of Lemma 6.8 in \cite{M}:

\begin{lemma}\label{6.6}
Let $S$ be a special subgroup of $H$. Let $g$ $\in$ $G\setminus H$. Let $h$ $\in$ $G\setminus g^S$. There exists a normal subgroup $L$ of $p$-power index in $H$ such that, if $\varphi$ : $H$ $\rightarrow$ $P$ = $H/L$ denotes the canonical projection, if $Q$ denotes the HNN extension of $P$ relative to $\varphi(K)$ and if $\overline{\varphi}$ : $G$ $\rightarrow$ $Q$ denotes the homomorphism induced by $\varphi$, we have $\overline{\varphi}(h)$ $\notin$ $\overline{\varphi}(g)^{\overline{\varphi}(S)}$.
\end{lemma}

\textit{Proof}: Write $g$ = $x_0t^{a_1}x_1\cdots t^{a_n}x_n$ and $h$ = $y_0t^{b_1}y_1\cdots t^{b_m}y_m$ in reduced forms. We have $n$ $\geq$ 1 -- as $g$ $\notin$ $H$.
\\

\underline{Step 1}: We assume that the first condition in Minasyan's criterion (see Lemma \ref{6.5}) is not satisfied by $g$ and $h$.
\\

It follows from Lemma 2.10 in \cite{Lo}, and Theorem \ref{6.1} (see, alternatively, \cite{Lo}, Theorem 2.11) that the special subgroup $K$ is closed in the pro-$p$ topology on $H$. (Note that this can also be obtained by combining Corollary \ref{5.11} and Theorem \ref{6.1}.) Thus, there exists a normal subgroup $L$ of $p$-power index in $H$ such that:
\begin{center}
$\forall$ $i$ $\in$ \{1,...,$n-1$\}, $x_i$ $\notin$ $KL$ ($\ast$), \\
$\forall$ $j$ $\in$ \{1,...,$m-1$\}, $y_j$ $\notin$ $KL$
($\ast\ast$).
\end{center}
We denote by $\varphi$ : $H$ $\rightarrow$ $P$ = $H/L$ the canonical projection. If $Q$ denotes the HNN extension of $P$ relative to $\varphi(K)$:
\begin{center}
$Q$ = $\langle$ $P$, $\overline{t}$ $\mid$ $\overline{t}^{-1}\varphi(k)\overline{t}$ = $\varphi(k)$, $\forall$ $k$ $\in$ $K$ $\rangle$,
\end{center}
and if $\overline{\varphi}$ : $G$ $\rightarrow$ $Q$ denotes the homomorphism induced by $\varphi$ -- with $\overline{\varphi}_{\mid_H}$ = $\varphi$ and $\overline{\varphi}(t)$ = $\overline{t}$ --, then $\overline{\varphi}(g)$ = $\overline{x_0}\overline{t}^{a_1}\overline{x_1}\cdots\overline{t}^{a_n}\overline{x_n}$ and $\overline{\varphi}(h)$ = $\overline{y_0}\overline{t}^{b_1}\overline{y_1}\cdots\overline{t}^{b_m}\overline{y_m}$ are reduced forms in $Q$ by ($\ast$) and ($\ast\ast$) -- where $\overline{x_i}$ = $\overline{\varphi}(x_i)$ ($i$ $\in$ \{0,...,$n$\}) and $\overline{y_j}$ = $\overline{\varphi}(y_j)$ ($j$ $\in$ \{0,...,$m$\}). But then the first condition in Minasyan's criterion will not hold for $\overline{\varphi}(g)$ and $\overline{\varphi}(h)$.
\\

\underline{Conclusion of Step 1}: We can assume that $m$ = $n$ and $a_i$ = $b_i$ for all $i$ $\in$ \{1,...,$n$\}.
\\

\underline{Step 2}: We assume that the second condition in Minasyan's criterion is not satisfied by $g$ and $h$. We set $x$ = $x_0\cdots x_n$ and $y$ = $y_0\cdots y_n$. Thus $y$ $\notin$ $x^S$.
\\

By the assumptions, $x^S$ is finitely $p$-separable in $H$. Therefore there exists a homomorphism $\varphi$ from $H$ onto a finite $p$-group $P$ such that $\varphi(y)$ $\notin$ $\varphi(x)^{\varphi(S)}$. Denote by $Q$ the HNN extension of $P$ relative to $\varphi(K)$, and by $\overline{\varphi}$ : $G$ $\rightarrow$ $Q$ the homomorphism induced by $\varphi$. Now let $f$ : $Q$ $\rightarrow$ $P$ be the natural homomorphism. We have:
\begin{center}
$f(\overline{\varphi}(g))$ = $f(\overline{x_0}\overline{t}^{a_1}\overline{x_1}\cdots\overline{t}^{a_n}\overline{x_n})$ = $\overline{x_0}\cdots\overline{x_n}$ = $\varphi(x)$, \\
$f(\overline{\varphi}(h))$ = $f(\overline{y_0}\overline{t}^{a_1}\overline{y_1}\cdots\overline{t}^{a_n}\overline{y_n})$ = $\overline{y_0}\cdots\overline{y_n}$ = $\varphi(y)$
\end{center}
(with the same notations as above). Since $\varphi(y)$ $\notin$ $\varphi(x)^{\varphi(S)}$, we see that $\overline{\varphi}(h)$ $\notin$ $\overline{\varphi}(g)^{\overline{\varphi}(S)}$.
\\

\underline{Conclusion of Step 2}: We can assume that $y$ $\in$ $x^S$. There exists $\alpha$ $\in$ $S$ such that $y$ = $\alpha x\alpha^{-1}$.
\\

\underline{End of the proof}: Considering Minasyan's criterion, since $h$ $\notin$ $g^S$, we must have:
\begin{center}
$\alpha C_S(x_0\cdots x_n)\cap y_0Kx_0^{-1}\cap (y_0y_1)K(x_0x_1)^{-1}\cap\cdots\cap(y_0\cdots y_{n-1})K(x_0\cdots x_{n-1})^{-1}$ = $\emptyset$.
\end{center}
As we noted above, $K$ is closed in the pro-$p$ topology on $H$; thus, there exists a normal subgroup $L$ of $p$-power index in $H$ such that:
\begin{center}
$\forall$ $i$ $\in$ \{1,...,$n-1$\}, $x_i$ $\notin$ $KL$ ($\ast$), \\
$\forall$ $j$ $\in$ \{1,...,$n-1$\}, $y_j$ $\notin$ $KL$ ($\ast\ast$),
\end{center}
Now by Lemma \ref{6.3}, there exists a normal subgroup $N$ of $p$-power index in $H$ such that $N$ $<$ $L$ and, if $\varphi$ : $H$ $\rightarrow$ $P$ = $H/N$ denotes the canonical projection, then:
\begin{center}
$\overline{\alpha}C_{\overline{S}}(\overline{x})\cap\overline{y_0}\overline{K}\overline{x_0}^{-1}\cap\overline{y_0}\hspace{1mm}\overline{y_1}\overline{K}(\overline{x_0}\hspace{1mm}\overline{x_1})^{-1}\cap\cdots\cap\overline{y_0}\cdots\overline{y_{n-1}}\overline{K}(\overline{x_0}\cdots\overline{x_{n-1}}) ^{-1}$ $\subset$ $\varphi((\alpha C_S(x)\cap y_0Kx_0^{-1}\cap y_0y_1K(x_0x_1)^{-1}\cap\cdots\cap y_0\cdots y_{n-1}K(x_0\cdots x_{n-1})^{-1})L)$ = $\emptyset$ ($\ast\ast\ast$).
\end{center}
where $\overline{S}$ = $\varphi(S)$, $\overline{\alpha}$ = $\varphi(\alpha)$, $\overline{x}$ = $\varphi(x)$, $\overline{x_i}$ = $\varphi(x_i)$ ($i$ $\in$ \{0,...,$n$\}), $\overline{y_j}$ = $\varphi(y_j)$ ($j$ $\in$ \{0,...,$n$\}). Let $Q$ be the HNN extension of $P$ relative to $\varphi(K)$ and let $\overline{\varphi}$ : $G$ $\rightarrow$ $Q$ be the homomorphism induced by $\varphi$. Then, by ($\ast$) and ($\ast\ast$), $\overline{\varphi}(g)$ = $\overline{x_0}\overline{t}^{a_1}\overline{x_1}\cdots\overline{t}^{a_n}\overline{x_n}$ and $\overline{\varphi}(h)$ = $\overline{y_0}\overline{t}^{a_1}\overline{y_1}\cdots\overline{t}^{a_n}\overline{y_n}$ are reduced forms in $Q$. So, in view of ($\ast\ast\ast$), we have $\overline{\varphi}(h)$ $\notin$ $\overline{\varphi}(g)^{\overline{\varphi}(S)}$.\hfill$\square$
\\

The following is the analogue of Lemma 8.8 in \cite{M}:

\begin{lemma}\label{6.7}
Let $g_0$ = $t^{a_1}x_1\cdots t^{a_n}x_n$ ($n$ $\geq$ 1) and $h_0$ = $t^{b_1}y_1\cdots t^{b_m}y_m$ be cyclically reduced elements of $G$. Let $h_1$,...,$h_k$ be elements of $G$. If $h_i$ $\notin$ $g_0^K$ for all $i$ $\in$ \{1,...,$k$\}, then there exists a normal subgroup $L$ of $p$-power index in $H$ such that, if $\varphi$ : $H$ $\rightarrow$ $P$ = $H/L$ denotes the canonical projection, if $Q$ denotes the HNN extension of $P$ relative to $\varphi(K)$, and if $\overline{\varphi}$ : $G$ $\rightarrow$ $Q$ denotes the homomorphism induced by $\varphi$, we have:
\begin{enumerate}
\item $\overline{\varphi}(g_0)$ = $\overline{t}^{a_1}\overline{x_1}\cdots\overline{t}^{a_n}\overline{x_n}$ and $\overline{\varphi}(h_0)$ = $\overline{t}^{b_1}\overline{y_1}\cdots\overline{t}^{b_m}\overline{y_m}$ are cyclically reduced in $Q$ -- where $\overline{x_i}$ = $\overline{\varphi}(x_i)$ ($i$ $\in$ \{1,...,$n$\}) and $\overline{y_j}$ = $\overline{\varphi}(y_j)$ ($j$ $\in$ \{1,...,$n$\}),
\item $\overline{\varphi}(h_i)$ $\notin$ $\overline{\varphi}(g_0)^{\overline{\varphi}(K)}$ for all $i$ $\in$ \{1,...,$k$\}.
\end{enumerate}
\end{lemma}

\textit{Proof}: As we noted above, $K$ is closed in the pro-$p$ topology on $H$; thus, there exists a normal subgroup $L_0$ of $p$-power index in $H$ such that:
\begin{center}
$\forall$ $i$ $\in$ \{1,...,$n-1$\}, $x_i$ $\notin$ $KL_0$ ($\ast$), \\
$\forall$ $j$ $\in$ \{1,...,$m-1$\}, $y_j$ $\notin$ $KL_0$ ($\ast\ast$).
\end{center}
Let $i$ $\in$ \{1,...,$k$\}. Since $h_i$ $\notin$ $g_0^K$, there exists a normal subgroup $L_i$ of $p$-power index in $H$ such that, if $\varphi_i$ : $H$ $\rightarrow$ $P_i$ = $H/L_i$ denotes the canonical projection, if $Q_i$ denotes the HNN extension of $P_i$ relative to $\varphi_i(K)$ and if $\overline{\varphi_i}$ : $G$ $\rightarrow$ $Q_i$ denotes the homomorphism induced by $\varphi_i$, we have $\overline{\varphi_i}(h_i)$ $\notin$ $\overline{\varphi_i}(g_0)^{\overline{\varphi_i}(K)}$ -- by Lemma \ref{6.6}. Set $L$ = $L_0\cap L_1\cdots\cap L_k$. Let $\varphi$ : $H$ $\rightarrow$ $P$ = $H/L$ be the canonical projection, let $Q$ be the HNN extension of $P$ relative to $\varphi(K)$, and let $\overline{\varphi}$ : $G$ $\rightarrow$ $Q$ be the homomorphism induced by $\varphi$. Since $L$ $<$ $L_0$, $\overline{\varphi}(g_0)$ = $\overline{t}^{a_1}\overline{x_1}\cdots\overline{t}^{a_n}\overline{x_n}$ and $\overline{\varphi}(h_0)$ = $\overline{t}^{b_1}\overline{y_1}\cdots\overline{t}^{b_m}\overline{y_m}$ are cyclically reduced in $Q$ by ($\ast$) and ($\ast\ast$) (with the same notations as in the statement of the lemma). As $L$ $<$ $L_i$ for all $i$ $\in$ \{1,...,$k$\}, we have $\overline{\varphi}(h_i)$ $\notin$ $\overline{\varphi}(g_0)^{\overline{\varphi}(K)}$ for all $i$ $\in$ \{1,...,$k$\}.\hfill$\square$

\begin{lemma}\label{6.8}
Let $G$ = $\langle$ $H$, $t$ $\mid$ $t^{-1}kt$ = $k$, $\forall$ $k$ $\in$ $K$ $\rangle$ be an HNN extension. Let $S$ be a subgroup of $H$. Let $g$ = $x_0t^{a_1}x_1\cdots t^{a_n}x_n$ be an element of $G$ in reduced form ($n$ $\geq$ 1). Then:
\begin{center}
$C_S(g)$ = $C_S(x_0\cdots x_n)\cap x_0Kx_0^{-1}\cap (x_0x_1)K(x_0x_1)^{-1}\cap\cdots\cap(x_0\cdots x_{n-1})K(x_0\cdots x_{n-1})^{-1}$.
\end{center}
\end{lemma}

\textit{Proof}: Proved in \cite{M} (see Lemma 7.12).\hfill$\square$
\\

The following is the analogue of Lemma 8.9 in \cite{M}:

\begin{lemma}\label{6.9}
Let $S$ be a special subgroup of $H$. Let $L$ be a normal subgroup of $p$-power index in $G$, and let $g$ = $x_0t^{a_1}x_1\cdots t^{a_n}x_n$ be an element of $G$ in reduced form and not contained in $H$. Then there exists a normal subgroup $N$ of $p$-power index of $H$ such that, if $\varphi$ : $H$ $\rightarrow$ $P$ = $H/N$ denotes the canonical projection, if $Q$ denotes the HNN extension of $P$ relative to $\varphi(K)$, and if $\overline{\varphi}$ : $G$ $\rightarrow$ $Q$ denotes the homomorphism induced by $\varphi$, we have:
\begin{enumerate}
\item $C_{\overline{\varphi}(S)}(\overline{\varphi}(g))$ $\subset$ $\overline{\varphi}(C_S(g)L)$,
\item $ker(\varphi)$ = $N$ $<$ $H\cap L$,
\item $ker(\overline{\varphi})$ $<$ $L$.
\end{enumerate}
\end{lemma}

\textit{Proof}: We have $n$ $\geq$ 1 -- as $g$ $\notin$ $H$.

As we noted above, $K$ is closed in the pro-$p$ topology on $H$. Therefore there exists a normal subgroup $M$ of $p$-power index in $H$ such that:
\begin{center}
$\forall$ $i$ $\in$ \{1,...,$n-1$\}, $x_i$ $\notin$ $KM$ ($\ast$).
\end{center}
We set $L'$ = $H\cap L$. Note that $L'$ is a normal subgroup of $p$-power index in $H$. Thus, up to replacing $M$ by $M\cap L'$, we can assume that $M$ $<$ $L'$. We set $x$ = $x_0\cdots x_n$. We have:
\begin{center}
$C_S(g)$ = $C_S(x)\cap x_0Kx_0^{-1}\cap(x_0x_1)K(x_0x_1)^{-1}\cap\cdots\cap(x_0\cdots x_{n-1})K(x_0\cdots x_{n-1})^{-1}$,
\end{center}
by Lemma \ref{6.8}. We denote by $I$ the intersection in the right-hand side. By Lemma \ref{6.3}, there exists a normal subgroup $N$ of $p$-power index in $H$ such that $N$ $<$ $M$ and, if $\varphi$ : $H$ $\rightarrow$ $P$ = $H/N$ denotes the canonical projection, we have:
\begin{center}
$C_{\overline{S}}(\overline{x})\cap\overline{x_0}\overline{K}\overline{x_0}^{-1}\cap\overline{x_0}\hspace{1mm}\overline{x_1}\overline{K}(\overline{x_0} \hspace{1mm}\overline{x_1})^{-1}\cap\cdots\cap\overline{x_0}\cdots\overline{x_{n-1}}\overline{K}(\overline{x_0}\cdots\overline{x_{n-1}})^{-1}$ $\subset$ $\varphi(IM)$
\end{center}
where $\overline{S}$ = $\varphi(S)$, $\overline{x}$ = $\varphi(x)$, $\overline{x_i}$ = $\varphi(x_i)$ ($i$ $\in$ \{0,...,$n-1$\}). We denote by $J$ the intersection in the left-hand side. Let $Q$ be the HNN extension of $P$ relative to $\varphi(K)$, and let $\overline{\varphi}$ : $G$ $\rightarrow$ $Q$ be the homomorphism induced by $\varphi$. Then $\overline{x_0}\overline{t}^{a_1}\overline{x_1}\cdots\overline{t}^{a_n}\overline{x_n}$ is a reduced form of $\overline{\varphi}(g)$ in $Q$ by ($\ast$). But then $C_{\overline{\varphi}(S)}(\overline{\varphi}(g))$ = $J$ -- by Lemma \ref{6.8}. Now $\varphi(M)$ $<$ $\varphi(L')$ = $\overline{\varphi}(L')$ $<$ $\overline{\varphi}(L)$. Therefore:
\begin{center}
$C_{\overline{\varphi}(S)}(\overline{\varphi}(g))$ = $J$ $\subset$ $\varphi(IM)$ = $\varphi(I)\varphi(M)$ $\subset$ $\overline{\varphi}(I)\overline{\varphi}(L)$ = $\overline{\varphi}(C_S(g))\overline{\varphi}(L)$ = $\overline{\varphi}(C_S(g)L)$.
\end{center}
Finally we remark that $ker(\varphi)$ = $N$ $<$ $M$ $<$ $L'$ = $H\cap L$ $<$ $L$. Since $ker(\overline{\varphi})$ is the normal closure of $ker(\varphi)$ in $G$, we conclude that $ker(\overline{\varphi})$ $<$ $L$ (because $L$ is normal in $G$).\hfill$\square$
\\

A \textit{prefix} of $t^{a_1}x_1\cdots t^{a_n}x_n$ is an element of $G$ of the form $t^{a_1}x_1\cdots t^{a_k}x_k$ for some $k$ $\in$ \{0,...,$n$\}. We need the following result:

\begin{proposition}\label{6.10}
Let $G$ = $\langle$ $H$, $t$ $\mid$ $t^{-1}kt$ = $k$, $\forall$ $k$ $\in$ $K$ $\rangle$ be an HNN extension. Let $g$ = $t^{a_1}x_1\cdots t^{a_n}x_n$ be a cyclically reduced element of $G$ ($n$ $\geq$ 1). Let \{$p_1$,...,$p_{n+1}$\} be the set of all prefixes of $g$ -- we are not assuming that $p_1$,...,$p_{n+1}$ are ordered. There are two cases:
\begin{enumerate}
\item if $x_n$ $\in$ $K$, then $n$ = 1 and $C_G(g)$ = $\langle t \rangle C_K(g)$;
\item if $x_n$ $\in$ $H\setminus K$, let \{$p_1$,...,$p_m$\} be the set of prefixes of $g$ satisfying $p_i^{-1}gp_i$ $\in$ $g^K$ ($m$ $\in$ \{0,...,$n$+1\}). For each $i$ $\in$ \{1,...,$m$\}, we choose $\alpha_i$ $\in$ $K$ such that $p_i^{-1}gp_i$ = $\alpha_i^{-1}g\alpha_i$. We set $\Omega$ = \{$\alpha_ip_i^{-1}$ $\mid$ $i$ $\in$ \{1,...,$m$\}\}. Then $C_G(g)$ = $C_K(g)\langle g\rangle\Omega$.
\end{enumerate}
\end{proposition}

\textit{Proof}: Proved in \cite{M} (see Proposition 7.8).\hfill$\square$
\\

The following is the analogue of Lemma 8.10 in \cite{M}:

\begin{lemma}\label{6.11}
Let $L$ be a normal subgroup of $p$-power index in $G$. Let $g_0$ = $t^{a_1}x_1\cdots t^{a_n}x_n$ ($n$ $\geq$ 1) be a cyclically reduced element of $G$. There exists a normal subgroup $N$ of $p$-power index in $H$ such that, if $\varphi$ : $H$ $\rightarrow$ $P$ = $H/N$ denotes the canonical projection, if $Q$ denotes the HNN extension of $P$ relative to $\varphi(K)$, and if $\overline{\varphi}$ : $G$ $\rightarrow$ $Q$ denotes the homomorphism induced by $\varphi$, we have:
\begin{enumerate}
\item $C_Q(\overline{\varphi}(g_0))$ $\subset$ $\overline{\varphi}(C_G(g_0)L)$,
\item $ker(\varphi)$ = $N$ $<$ $H\cap L$,
\item $ker(\overline{\varphi})$ $<$ $L$.
\end{enumerate}
\end{lemma}

\textit{Proof}: Let \{$p_1$,...,$p_{n+1}$\} be the set of all prefixes of $g_0$. Renumbering $p_1$,...,$p_{n+1}$, if necessary, we can assume that there exists $m$ $\in$ \{1,...,$n+1$\} such that $p_i^{-1}g_0p_i$ $\in$ $g_0^K$ for all $i$ $\in$ \{1,...,$m$\}, and $p_i^{-1}g_0p_i$ $\notin$ $g_0^K$ for all $i$ $\in$ \{$m$+1,...,$n$+1\}. For each $i$ $\in$ \{1,...,$m$\}, we choose $\alpha_i$ $\in$ $K$ such that $p_i^{-1}g_0p_i$ = $\alpha_i^{-1}g_0\alpha_i$. We set $\Omega$ = \{$\alpha_ip_i^{-1}$ $\mid$ $i$ $\in$ \{1,...,$m$\}\}. We set $h_i$ = $p_i^{-1}g_0p_i$ for all $i$ $\in$ \{$m$+1,...,$n$+1\}. By Lemma \ref{6.7}, there exists a normal subgroup $N_1$ of $p$-power index in $H$ such that, if $\varphi_1$ : $H$ $\rightarrow$ $P_1$ = $H/N_1$ denotes the canonical projection, if $Q_1$ denotes the HNN extension of $P_1$ relative to $\varphi_1(K)$, and if $\overline{\varphi_1}$ : $G$ $\rightarrow$ $Q_1$ denotes the homomorphism induced by $\varphi_1$, then $\varphi_1(g_0)$ is cyclically reduced in $Q_1$, and $\overline{\varphi_1}(h_i)$ $\notin$ $\overline{\varphi_1}(g_0)^{\overline{\varphi_1}(K)}$ for all $i$ $\in$ \{$m$+1,...,$n$+1\}. On the other hand, by Lemma \ref{6.9}, there exists a normal subgroup $N_2$ of $p$-power index in $H$ such that, if $\varphi_2$ : $H$ $\rightarrow$ $P_2$ = $H/N_2$ denotes the canonical projection, if $Q_2$ denotes the HNN extension of $P_2$ relative to $\varphi_2(K)$ and if $\overline{\varphi_2}$ : $G$ $\rightarrow$ $Q_2$ denotes the homomorphism induced by $\varphi_2$, we have: $C_{\overline{\varphi_2}(K)}(\overline{\varphi_2}(g_0))$ $\subset$ $\overline{\varphi_2}(C_K(g_0)L)$, $ker(\varphi_2)$ $<$ $H\cap L$, and $ker(\overline{\varphi_2})$ $<$ $L$. Set $N$ = $N_1\cap N_2$. Let $\varphi$ : $H$ $\rightarrow$ $P$ = $H/N$ be the canonical projection, let $Q$ be the HNN extension of $P$ relative to $\varphi(K)$, and let $\overline{\varphi}$ : $G$ $\rightarrow$ $Q$ be the homomorphism induced by $\varphi$. Since $N$ $<$ $N_1$, $\overline{\varphi}(g_0)$ is cyclically reduced in $Q$ and $\overline{\varphi}(h_i)$ $\notin$ $\overline{\varphi}(g_0)^{\overline{\varphi}(K)}$ for all $i$ $\in$ \{$m$+1,...,$n$+1\}. On the other hand, since $N$ $<$ $N_2$, we have:
\begin{center}
$\overline{\varphi}^{-1}(C_{\overline{\varphi}(K)}(\overline{\varphi}(g_0)))$ $\subset$ $\overline{\varphi_2}^{-1}(C_{\overline{\varphi_2}(K)}(\overline{\varphi_2}(g_0)))$ $\subset$ $C_K(g_0)L$ ($\ast$).
\end{center}
There are two cases:
\\

\underline{Case 1}: $x_n$ $\in$ $K$. Then $n$ = 1, $C_G(g_0)$ = $\langle t \rangle C_K(g_0)$, and $C_Q(\overline{\varphi}(g_0))$ = $\langle\overline{t}\rangle C_{\varphi(K)}(\overline{\varphi}(g_0))$ -- by Proposition \ref{6.10}. Now ($\ast$) implies:
\begin{center}
$C_Q(\overline{\varphi}(g_0))$ $\subset$ $\langle\overline{\varphi}(t)\rangle\overline{\varphi}(C_K(g_0)L)$ = $\overline{\varphi}(\langle t\rangle C_K(g_0)L)$ = $\overline{\varphi}(C_G(g_0)L)$.
\end{center}

\vspace{2mm}

\underline{Case 2}: $x_n$ $\in$ $H\setminus K$. If $i$ $\in$ \{1,...,$m$\}, $\overline{\varphi}(p_i)^{-1}\overline{\varphi}(g_0)\overline{\varphi}(p_i)$ = $\overline{\varphi}(p_i^{-1}g_0p_i)$ $\in$ $\overline{\varphi}(g_0)^{\varphi(K)}$ -- because $p_i^{-1}g_0p_i$ $\in$ $g_0^K$ --, whereas if $i$ $\in$ \{$m$+1,...,$n$+1\}, $\overline{\varphi}(p_i)^{-1}\\ \overline{\varphi}(g_0)$$\overline{\varphi}(p_i)$ = $\overline{\varphi}(h_i)$ $\notin$ $\overline{\varphi}(g_0)^{\varphi(K)}$. Therefore \{$\overline{\varphi}(p_1)$,...,$\overline{\varphi}(p_m)$\} is the set of all prefixes of $\overline{\varphi}(g_0)$ satisfying $\overline{\varphi}(p_i)^{-1}\overline{\varphi}(g_0)\overline{\varphi}(p_i)$ $\in$ $\overline{\varphi}(g_0)^{\varphi(K)}$. By Proposition \ref{6.10}, $C_G(g_0)$ = $C_K(g_0)\langle g_0 \rangle\Omega$, and $C_Q(\overline{\varphi}(g_0))$ = $C_{\varphi(K)}(\overline{\varphi}(g_0))\langle\overline{\varphi}(g_0)\rangle\overline{\Omega}$, where $\overline{\Omega}$ = $\overline{\varphi}(\Omega)$ = \{$\overline{\varphi}(\alpha_i)\overline{\varphi}(p_i)^{-1}$ $\mid$ $i$ $\in$ \{1,...,$m$\}\}. We deduce that:
\begin{center}
$C_Q(\overline{\varphi}(g_0))$ $\subset$ $\overline{\varphi}(C_K(g_0)L)\langle\overline{\varphi}(g_0)\rangle\overline{\varphi}(\Omega)$ = $\overline{\varphi}(C_K(g_0)L\langle g_0\rangle\Omega)$ = $\overline{\varphi}(C_G(g_0)L)$.
\end{center}
\hfill$\square$

\begin{proposition}\label{6.12}
Let $G$ be a right-angled Artin group of rank $r$ ($r$ $\geq$ 1). Let $g$ $\in$ $G$. If $g$ $\neq$ 1, there exists a special subgroup $H$ of rank $r-1$ of $G$ such that $g$ $\notin$ $H^G$, where $H^G$ = $\cup_{h\in H}h^G$.
\end{proposition}

\textit{Proof}: Proved in \cite{M} (see Lemma 6.8).\hfill$\square$

\begin{lemma}\label{6.13}
Every special subgroup $S$ of $G$ satisfies the $p$-centralizer condition in $G$ $(pCC_G)$.
\end{lemma}

\textit{Proof}: Let $g$ $\in$ $G$. Let $L$ be a normal subgroup of $p$-power index in $G$. There are two cases: 
\\

\underline{Case 1}: $S$ $\neq$ $G$.

Let $H$ be a special subgroup of rank $r-1$ of $G$ such that $S$ $<$ $H$. Then $G$ can be written as an HNN extension of $H$, relative to a special subgroup $K$ of $H$:
\begin{center}
$G$ = $\langle$ $H$, $t$ $\mid$ $t^{-1}kt$ = $k$, $\forall$ $k$ $\in$ $K$ $\rangle$.
\end{center}
We set $L'$ = $H\cap L$. We note that $L'$ is a normal subgroup of $p$-power index in $H$. There are two cases:

\vspace{1mm}

\underline{Subcase 1}: $g$ $\in$ $H$. By the assumptions, the pair ($S$,$g$) satisfies the $p$-centralizer condition in $H$ ($pCC_H$). There exists a normal subgroup $M$ of $p$-power index in $H$ such that $M$ $<$ $L'$ and, if $\psi$ : $H$ $\rightarrow$ $P$ = $H/M$ denotes the canonical projection, we have:
\begin{center}
$C_{\psi(S)}(\psi(g))$ $\subset$ $\psi(C_S(g)L')$ ($\ast$).
\end{center}
We denote by $f$ : $G$ $\rightarrow$ $H$ the natural homomorphism. We note that $f^{-1}(M)$ is a normal subgroup of $p$-power index in $G$ (because $f^{-1}(M)$ is the kernel of the homomorphism $\psi\circ f$). Therefore, $N$ = $L\cap f^{-1}(M)$ is a normal subgroup of $p$-power index in $G$. Moreover $N$ $<$ $L$ and $f(N)$ $<$ $M$. We denote by $\varphi$ : $G$ $\rightarrow$ $Q$ = $G/N$ the canonical projection. We observe that $ker(\psi)$ = $M$, $ker(\varphi)$ = $N$, $M$ $<$ $f^{-1}(M)\cap L\cap H$ = $N\cap H$, and $N\cap H$ $\subset$ $f(N)$ $<$ $M$. Therefore $M$ = $N\cap H$. Thus we can assume that $P$ $<$ $Q$ and $\varphi_{\mid_H}$ = $\psi$. But then $\psi(L')$ = $\varphi(L')$ $\subset$ $\varphi(L)$. Recall that $g$ $\in$ $H$ and $S$ $<$ $H$. Thus considering ($\ast$), we obtain:
\begin{center}
$C_{\varphi(S)}(\varphi(g))$ = $C_{\psi(S)}(\psi(g))$ $\subset$ $\psi(C_S(g))\psi(L')$ $\subset$ $\varphi(C_S(g))\varphi(L)$ = $\varphi(C_S(g)L)$.
\end{center}

\vspace{1mm}

\underline{Subcase 2}: $g$ $\in$ $G\setminus H$. Write $g$ = $x_0t^{a_1}x_1\cdots t^{a_n}x_n$ in a reduced form ($n$ $\geq$ 1). Then, by Lemma \ref{6.9}, there exists a normal subgroup $M$ of $p$-power index in $H$ such that, if $\psi$ : $H$ $\rightarrow$ $P$ = $H/M$ denotes the canonical projection, if $Q$ denotes the HNN extension of $P$ relative to $\psi(K)$, and if $\overline{\psi}$ : $G$ $\rightarrow$ $Q$ denotes the homomorphism induced by $\psi$, then: $C_{\overline{\psi}(S)}(\overline{\psi}(g))$ $\subset$ $\overline{\psi}(C_S(g)L)$, $ker(\psi)$ $<$ $H\cap L$, and $ker(\overline{\psi})$ $<$ $L$. We note that $\overline{\psi}(S)\cap\overline{\psi}(L)$ = $\psi(S)\cap\overline{\psi}(L)$ $<$ $P$ is finite. Since $Q$ is residually $p$-finite, $\overline{\psi}(g)^{\overline{\psi}(S)\cap\overline{\psi}(L)}$ is finitely $p$-separable in $Q$. Therefore, by Lemma \ref{3.5}, there exists a normal subgroup $N$ of $p$-power index in $Q$ such that $N$ $<$ $\overline{\psi}(L)$ and, if $\chi$ : $Q$ $\rightarrow$ $R$ = $Q/N$ denotes the canonical projection, then:
\begin{center}
$C_{\chi(\overline{\psi}(S))}(\chi(\overline{\psi}(g)))$ $\subset$ $\chi(C_{\overline{\psi}(S)}(\overline{\psi}(g))\overline{\psi}(L))$.
\end{center}
We set $\varphi$ = $\chi\circ\overline{\psi}$ : $G$ $\rightarrow$ $R$. We have: $ker(\varphi)$ = $\overline{\psi}^{-1}(ker(\chi))$ = $\overline{\psi}^{-1}(N)$ $\subset$ $\overline{\psi}^{-1}(\overline{\psi}(L))$ = $Lker(\overline{\psi})$. Now $ker(\overline{\psi})$ $<$ $L$. Then $ker(\varphi)$ $<$ $L$. And:
\begin{center}
$C_{\varphi(S)}(\varphi(g))$ = $C_{\chi(\overline{\psi}(S))}(\chi(\overline{\psi}(g)))$ $\subset$ $\chi(C_{\overline{\psi}(S)}(\overline{\psi}(g))\overline{\psi}(L))$ $\subset$ $\chi(\overline{\psi}(C_S(g)L)\overline{\psi}(L))$ = $\varphi(C_S(g)L)$.
\end{center}

\vspace{2mm}

\underline{Case 2}: $S$ = $G$.

If $g$ = 1, then the result is trivial. Thus we can assume that $g$ $\neq$ 1. Then, by Proposition \ref{6.12}, there exists a special subgroup $H$ of rank $r-1$ of $G$ such that $g$ $\notin$ $H^G$. As above, $G$ can be written as an HNN extension of $H$ relative to a special subgroup $K$ of $H$:
\begin{center}
$G$ = $\langle$ $H$, $t$ $\mid$ $t^{-1}kt$ = $k$, $\forall$ $k$ $\in$ $K$ $\rangle$.
\end{center}
Let $g_0$ = $t^{a_1}x_1\cdots t^{a_n}x_n$ be a cyclically reduced element in $G$ conjugate to $g$. Choose $\alpha$ $\in$ $G$ such that $g$ = $\alpha g_0 \alpha^{-1}$. Note that $g$ $\notin$ $H^G$ implies that $n$ $\geq$ 1. By Lemma \ref{6.11}, there exists a normal subgroup $M$ of $p$-power index in $H$ such that, if $\psi$ : $H$ $\rightarrow$ $P$ = $H/M$ denotes the canonical projection, if $Q$ denotes the HNN extension of $P$ relative to $\psi(K)$, and if $\overline{\psi}$ : $G$ $\rightarrow$ $Q$ denotes the homomorphism induced by $\psi$, then: $C_Q(\overline{\psi}(g_0))$ $\subset$ $\overline{\psi}(C_G(g_0)L)$, $ker(\psi)$ $<$ $H\cap L$, and $ker(\overline{\psi})$ $<$ $L$. Now $Q$ is hereditarily conjugacy $p$-separable by Corollary \ref{4.3}. Then $Q$ satisfies the $p$-centralizer condition by Proposition \ref{3.6}. There exists a normal subgroup $N$ of $p$-power index in $Q$ such that $N$ $<$ $\overline{\psi}(L)$, and if $\chi$ : $Q$ $\rightarrow$ $R$ = $Q/N$ denotes the canonical projection, we have:
\begin{center}
$C_R(\chi(\overline{\psi}(g_0)))$ $\subset$ $\chi(C_Q(\overline{\psi}(g_0))\overline{\psi}(L))$.
\end{center}
We set $\varphi$ = $\chi\circ\overline{\psi}$ : $G$ $\rightarrow$ $R$. As above, we have $ker(\varphi)$ = $\overline{\psi}^{-1}(ker(\chi))$ = $\overline{\psi}^{-1}(N)$ $\subset$ $\overline{\psi}^{-1}(\overline{\psi}(L))$ = $Lker(\overline{\psi})$. Now $ker(\overline{\psi})$ $<$ $L$. Then $ker(\varphi)$ $<$ $L$. And:
\begin{center}
$C_R(\varphi(g_0))$ = $C_{\varphi(G)}(\varphi(g_0))$ = $C_{\chi(\overline{\psi}(G))}(\chi(\overline{\psi}(g_0)))$ $\subset$ $\chi(C_{\overline{\psi}(G)}(\overline{\psi}(g_0))\overline{\psi}(L))$ $\subset$ $\chi(\overline{\psi}(C_G(g_0)L)\overline{\psi}(L))$ = $\varphi(C_G(g_0)L)$.
\end{center}
Finally:
\begin{center}
$\varphi(\alpha)C_R(\varphi(g_0))\varphi(\alpha)^{-1}$ $\subset$ $\varphi(\alpha)\varphi(C_G(g_0)L)\varphi(\alpha)^{-1}$.
\end{center}
That is,
\begin{center}
$C_R(\varphi(g))$ $\subset$ $\varphi(C_G(g)L)$.
\end{center}
\hfill$\square$

\begin{lemma}\label{6.14}
For every $g$ $\in$ $G$ and for every special subgroup $S$ of $G$, $g^S$ is finitely $p$-separable in $G$.
\end{lemma}

\textit{Proof}: There are two cases: 
\\

\underline{Case 1}: $S$ $\neq$ $G$.

Let $H$ be a special subgroup of rank $r-1$ of $G$ such that $S$ $<$ $H$. As above, $G$ can be written as an HNN extension of $H$ relative to a special subgroup $K$ of $H$:
\begin{center}
$G$ = $\langle$ $H$, $t$ $\mid$ $t^{-1}kt$ = $k$, $\forall$ $k$ $\in$ $K$ $\rangle$.
\end{center}
Let $g$ $\in$ $G$. There are two cases:

\vspace{1mm}

\underline{Subcase 1}: $g$ $\in$ $H$. Then $g^S$ is finitely $p$-separable in $H$ by the assumptions. Since $G$ is residually $p$-finite by Theorem \ref{6.1}, $g^S$ is finitely $p$-separable in $G$ by Lemma \ref{5.12}.

\vspace{1mm}

\underline{Subcase 2}: $g$ $\in$ $G\setminus H$. Let $h$ $\in$ $G\setminus g^S$. By Lemma \ref{6.6}, there exists a normal subgroup $L$ of $p$-power index in $H$ such that, if $\psi$ : $H$ $\rightarrow$ $P$ = $H/L$ denotes the canonical projection, if $Q$ denotes the HNN extension of $P$ relative to $\psi(K)$, and if $\overline{\psi}$ : $G$ $\rightarrow$ $Q$ denotes the homomorphism induced by $\psi$, we have $\overline{\psi}(h)$ $\notin$ $\overline{\psi}(g)^{\overline{\psi}(S)}$. Now $\overline{\psi}(S)$ = $\psi(S)$ $<$ $P$ is finite and $Q$ is residually $p$-finite. Then there exists a homomorphism $\chi$ : $Q$ $\rightarrow$ $R$ from $Q$ onto a finite $p$-group $R$ such that $\chi(\overline{\psi}(h))$ $\notin$ $\chi(\overline{\psi}(g)^{\overline{\psi}(S)})$. Thus the homomorphism $\varphi$ = $\chi\circ\overline{\psi}$ : $G$ $\rightarrow$ $R$ satisfies the condition $\varphi(h)$ $\notin$ $\varphi(g^S)$, as required. 
\\

\underline{Case 2}: $S$ = $G$. 

Let $g$ $\in$ $G$. If $g$ = 1, then, since $G$ is residually $p$-finite by Theorem \ref{6.1}, $g^G$ = \{1\} is finitely $p$-separable in $G$. Thus we can assume that $g$ $\neq$ 1. Then, by Proposition \ref{6.12}, there exists a special subgroup $H$ of rank $r-1$ of $G$ such that $g$ $\notin$ $H^G$. As above, $G$ can be written as an HNN extension of $H$ relative to a special subgroup $K$ of $H$:
\begin{center}
$G$ = $\langle$ $H$, $t$ $\mid$ $t^{-1}kt$ = $k$, $\forall$ $k$ $\in$ $K$ $\rangle$.
\end{center}
Let $h$ $\in$ $G\setminus g^G$. Let $g_0$ = $t^{a_1}x_1\cdots t^{a_n}x_n$ and $h_0$ = $t^{b_1}y_1\cdots t^{b_m}y_m$ be cyclically reduced elements of $G$ conjugate to $g$ and $h$ respectively. Note that $g$ $\notin$ $H^G$ implies that $n$ $\geq$ 1. There are two cases:

\vspace{1mm}

\underline{Subcase 1}: $h_0$ $\in$ $H$.  Then, by Lemma \ref{6.7}, there exists a normal subgroup $L$ of $p$-power index in $H$ such that, if $\psi$ : $H$ $\rightarrow$ $P$ = $H/L$ denotes the canonical projection, if $Q$ denotes the HNN extension of $P$ relative to $\psi(K)$, and if $\overline{\psi}$ : $G$ $\rightarrow$ $Q$ denotes the homomorphism induced by $\psi$, then: $\overline{\psi}(g_0)$ = $\overline{t}^{a_1}\overline{x_1}\cdots\overline{t}^{a_n}\overline{x_n}$ is cyclically reduced in $Q$ -- where $\overline{x_i}$ = $\overline{\psi}(x_i)$ ($i$ $\in$ \{1,...,$n$\}). Since $n$ $\geq$ 1, we have: $\overline{\psi}(g_0)$ $\notin$ $P^Q$ = $\overline{\psi}(H^G)$. Therefore $\overline{\psi}(g_0)$ $\notin$ $\overline{\psi}(h_0)^Q$ = $\overline{\psi}(h_0^G)$ $\subset$ $\overline{\psi}(H^G)$. Now $Q$ is conjugacy $p$-separable by Corollary \ref{4.3}. Then there exists a homomorphism $\chi$ from $Q$ onto a finite $p$-group $R$ such that $\chi(\overline{\psi}(g_0))$ $\notin$ $\chi(\overline{\psi}(h_0))^R$. Therefore $\chi(\overline{\psi}(g))$ $\notin$ $\chi(\overline{\psi}(h))^R$. Thus the homomorphism $\varphi$ = $\chi\circ\overline{\psi}$ : $G$ $\rightarrow$ $R$ satisfies the condition $\varphi(h)$ $\notin$ $\varphi(g^S)$, as desired.

\vspace{1mm}

\underline{Subcase 2}: $h_0$ $\in$ $G\setminus H$. Let \{$h_1$,...,$h_m$\} be the set of all cyclic permutations of $h_0$. Then, since $h$ $\notin$ $g^G$, we have: $h_i$ $\notin$ $g_0^G$ for all $i$ $\in$ \{1,...,$m$\}. Therefore, by Lemma \ref{6.7}, there exists a normal subgroup $L$ of $p$-power index in $H$ such that, if $\psi$ : $H$ $\rightarrow$ $P$ = $H/L$ denotes the canonical projection, if $Q$ denotes the HNN extension of $P$ relative to $\psi(K)$, and if $\overline{\psi}$ : $G$ $\rightarrow$ $Q$ denotes the homomorphism induced by $\psi$, then: $\overline{\psi}(g_0)$ = $\overline{t}^{a_1}\overline{x_1}\cdots\overline{t}^{a_n}\overline{x_n}$ and $\overline{\psi}(h_0)$ = $\overline{t}^{b_1}\overline{y_1}\cdots\overline{t}^{b_m}\overline{y_m}$ are cyclically reduced in $Q$ -- where $\overline{x_i}$ = $\overline{\psi}(x_i)$ ($i$ $\in$ \{1,...,$n$\}) and $\overline{y_j}$ = $\overline{\psi}(y_j)$ ($j$ $\in$ \{1,...,$n$\}) -- and $\overline{\psi}(h_i)$ $\notin$ $\overline{\psi}(g_0)^{\overline{\psi}(K)}$ for all $i$ $\in$ \{1,...,$m$\}. Consequently, by Lemma \ref{2.3}, $\overline{\psi}(g_0)$ $\notin$ $\overline{\psi}(h_0)^Q$. Now $Q$ is conjugacy $p$-separable by Corollary \ref{4.3}. Then there exists a homomorphism $\chi$ from $Q$ onto a finite $p$-group $R$ such that: $\chi(\overline{\psi}(g_0))$ $\notin$ $\chi(\overline{\psi}(h_0))^R$. Hence $\chi(\overline{\psi}(g))$ $\notin$ $\chi(\overline{\psi}(h))^R$. Thus the homomorphism $\varphi$ = $\chi\circ\overline{\psi}$ : $G$ $\rightarrow$ $R$ satisfies the condition $\varphi(h)$ $\notin$ $\varphi(g^S)$, as required.\hfill$\square$
\\

\textit{Proof of Proposition \ref{6.2}}: We argue by induction on the rank $r$ of $G$. If $r$ = 0, then the result is trivial. Thus we can assume that $r$ $\geq$ 1 and that the result has been proved for 1,...,$r-1$. Now, Proposition \hyperlink{6.2.1}{6.2.1} follows from Lemma \ref{6.13}, and Proposition \hyperlink{6.2.2}{6.2.2} follows from Lemma \ref{6.14}.\hfill$\square$
\\

We are now ready to prove:

\begin{theorem}\label{6.15}
Every right-angled Artin group is hereditarily conjugacy $p$-separable.
\end{theorem}

\textit{Proof}: Let $G$ be a right-angled Artin group. Let $g$ $\in$ $G$. Then $g^G$ is finitely $p$-separable in $G$ by Proposition \hyperlink{6.2.1}{6.2.1}. We deduce that $G$ is conjugacy $p$-separable. On the other hand, $G$ satisfies the $p$-centralizer condition by Proposition \hyperlink{6.2.2}{6.2.2}. We conclude that $G$ is hereditarily conjugacy $p$-separable by Proposition \ref{3.6}.\hfill$\square$

\section{Applications}

\hspace{5mm}The first application that we mention is an application of our main theorem to separability properties of $G_\Gamma$.
\\

For a group $G$, we denote by $(C^n(G))_{n \geq 1}$ the lower central series of $G$. Recall that $(C^n(G))_{n \geq 1}$ is defined inductively by $C^1(G)$ = $G$, and $C^{n+1}(G)$ = $[G,C^n(G)]$ for all $n$ $\geq$ 1.

\begin{corollary}\label{7.1}
Every right-angled Artin group is conjugacy separable in the class of torsion-free nilpotent groups.
\end{corollary}

\textit{Proof}: Let $G$ be a right-angled Artin group. Let $g$, $h$ $\in$ $G$ such that $g$ $\nsim$ $h$. Let $p$ be a prime number. Then $G$ is conjugacy $p$-separable by Theorem \ref{6.15}. Thus, there exists a homomorphism $\varphi$ from $G$ onto a finite $p$-group $P$ such that $\varphi(g)$ $\nsim$ $\varphi(h)$. Now, $P$ is nilpotent. Therefore, there exists $n$ $\geq$ 1 such that $C^n(P)$ = \{1\}. Let $\pi$ : $G$ $\rightarrow$ $\frac{G}{C^n(G)}$ be the canonical projection. It follows from \cite{DK2}, Theorem 2.1, that for all $n$ $\geq$ 1, there exists $d_n$ $\in$ $\mathbbm{N}$ such that:
\begin{center}
$\frac{C^n(G)}{C^{n+1}(G)}$ $\simeq$ $\mathbbm{Z}^{d_n}$.
\end{center}
Thus, an easy induction on $n$ shows that $\frac{G}{C^n(G)}$ is torsion-free for all $n$ $\geq$ 1. Hence $\frac{G}{C^n(G)}$ is a torsion-free nilpotent group for all $n$ $\geq$ 1. Since $\varphi(C^n(G))$ $<$ $C^n(P)$ = \{1\}, $\varphi$ induces a homomorphism $\widetilde{\varphi}$ : $\frac{G}{C^n(G)}$ $\rightarrow$ $P$ such that $\varphi$ = $\widetilde{\varphi}\circ\pi$. As $\varphi(g)$ $\nsim$ $\varphi(h)$, we have $\pi(g)$ $\nsim$ $\pi(h)$.\hfill$\square$
\\

We now turn to applications of our main theorem to residual properties of $Out(G_\Gamma)$.
\\

An automorphism $\varphi$ of a group $G$ is said to be \textit{conjugating} if for every $g$ $\in$ $G$, $\varphi$($g$) $\sim$ $g$. We say that $G$ has \textit{Property A} if every conjugating automorphism of $G$ is inner. The following proposition is due to Minasyan (see \cite{M}, Proposition 6.9):

\begin{proposition}\label{7.2}
Right-angled Artin groups have Property A.
\end{proposition}

For a group $G$, we denote by $\mathcal{I}_p(G)$ the kernel of the natural homomorphism $Out(G)$ $\rightarrow$ $GL(H_1(G,\mathbbm{F}_p))$ (where $\mathbbm{F}_p$ denotes the finite field with $p$ elements). The following theorem is due to Paris (see \cite{P}, Theorem 2.5):

\begin{theorem}\label{7.3}
Let $G$ be a finitely generated group. If $G$ is conjugacy $p$-separable and has Property A, then $\mathcal{I}_p(G)$ is residually $p$-finite.
\end{theorem}

Thus, combining Theorem \ref{7.3} and Proposition \ref{7.2} with Theorem \ref{6.15}, we obtain:

\begin{corollary}\label{7.4}
The outer automorphism group of a right-angled Artin group is virtually residually $p$-finite.
\end{corollary}

The following theorem is due to Myasnikov (see \cite{My}, Theorem 1):

\begin{theorem}\label{7.5}
Let $G$ be a finitely generated group. If $G$ is conjugacy $p$-separable and has property A, then $Out(G)$ is residually $\mathcal{K}$, where $\mathcal{K}$ is the class of all outer automorphism groups of finite $p$-groups.
\end{theorem}

Thus, combining Theorem \ref{7.5} and Proposition \ref{7.2} with Theorem \ref{6.15}, we obtain:

\begin{corollary}\label{7.6}
The outer automorphism group of a right-angled Artin group is residually $\mathcal{K}$, where $\mathcal{K}$ is the class of all outer automorphism groups of finite $p$-groups.
\end{corollary}

In the remainder of this paper, we prove Theorem \ref{7.14}. Let $G$ = $G_\Gamma$ be a right-angled Artin group. Let $r$ be the rank of $G$. We denote by $T(G)$ the kernel of the natural homomorphism $Aut(G)$ $\rightarrow$ $GL_r(\mathbbm{Z})$, and by $\mathcal{T}(G)$ the kernel of the natural homomorphism $Out(G)$ $\rightarrow$ $GL_r(\mathbbm{Z})$. Note that $\mathcal{T}(G)$ = $T(G)/Inn(G)$. Day proved that $T(G)$ is finitely generated (see \cite{D2}, Theorem B). Therefore $\mathcal{T}(G)$ is finitely generated. 
\\

In order to prove Theorem \ref{7.14}, we have to introduce the notion of separating $\mathbbm{Z}$-linear central filtration.

Recall that a \textit{central filtration} on a group $G$ is a sequence $(G_n)_{n \geq 1}$ of subgroups of $G$ satisfying the conditions:

\begin{center}
$G_1$ = $G$, \\
$G_n$ $>$ $G_{n+1}$, \\
$[G_m,G_n]$ $<$ $G_{m+n}$ for all $m$, $n$ $\geq$ 1.
\end{center}

Let $\mathcal{F}$ = $(G_n)_{n \geq 1}$ be a central filtration. Then the mapping $G\times G$ $\rightarrow$ $G$, $(x,y)$ $\mapsto$ $xyx^{-1}y^{-1}$ induces on:
\begin{center}
$\mathcal{L}_\mathcal{F}(G)$ = $\bigoplus_{n \geq 1}\frac{G_n}{G_{n+1}}$
\end{center}
a Lie bracket which makes $\mathcal{L}_\mathcal{F}(G)$ into a graded Lie $\mathbbm{Z}$-algebra.

We say that $(G_n)_{n \geq 1}$ is a \textit{separating filtration} if $\cap_{n \geq 1}G_n$ = \{1\}. We say that $(G_n)_{n \geq 1}$ is \textit{$\mathbbm{Z}$-linear} if for all $n$ $\geq$ 1, the $\mathbbm{Z}$-module $\frac{G_n}{G_{n+1}}$ is free of finite rank.
\\

For a group $G$, we denote by $(C^n_\mathbbm{Z}(G))_{n \geq 1}$ the sequence of subgroups of $G$ defined inductively by $C^1_\mathbbm{Z}(G)$ = $G$, $[G,C^n_\mathbbm{Z}(G)]$ $<$ $C^{n+1}_\mathbbm{Z}(G)$, and $\frac{C^{n+1}_\mathbbm{Z}(G)}{[G,C^n_\mathbbm{Z}(G)]}$ is the torsion subgroup of $\frac{C^n_\mathbbm{Z}(G)}{[G,C^n_\mathbbm{Z}(G)]}$ for all $n$ $\geq$ 1.

\begin{proposition}
For all $m$, $n$ $\geq$ 1, $[C^m_\mathbbm{Z}(G),C^n_\mathbbm{Z}(G)]$ $<$ $C^{m+n}_\mathbbm{Z}(G)$.
\end{proposition}

\textit{Proof}: Proved in \cite{BL} (see Proposition 7.2).\hfill$\square$
\\

Thus, $(C^n_\mathbbm{Z}(G))_{n \geq 1}$ is a central filtration on $G$. We denote by $\mathcal{L}_\mathbbm{Z}(G)$ the corresponding graded Lie $\mathbbm{Z}$-algebra.
\\

For a Lie algebra $\mathfrak{g}$, we denote by $Z(\mathfrak{g})$ the center of $\mathfrak{g}$. Let $G$ be a group. For $n$ $\geq$ 1, we denote by $A_n$ the kernel of the natural homomorphism $Aut(G)$ $\rightarrow$ $Aut(\frac{G}{C^{n+1}_\mathbbm{Z}(G)})$. Let $\pi$ : $Aut(G)$ $\rightarrow$ $Out(G)$ be the canonical projection. For $n$ $\geq$ 1, we set $B_n$ = $\pi(G_n)$.

\begin{theorem}\label{7.8}
If $G^{ab}$ is finitely generated, and $Z(\mathbbm{F}_p \otimes \mathcal{L}_\mathbbm{Z}(G))$ = \{0\} for every prime number $p$, then $(B_n)_{n \geq 1}$ is a $\mathbbm{Z}$-linear central filtration on $B_1$. Furthermore, $(B_n)_{n \geq 1}$ is separating if and only if $G$ satisfies the condition:
\begin{center}
$(IN(G))$: For every $\varphi$ $\in$ $Aut(G)$, if $\varphi$ induces an inner automorphism of $\frac{G}{C^n_\mathbbm{Z}(G)}$ for all $n$ $\geq$ 1, then $\varphi$ is inner.
\end{center}
\end{theorem}

\textit{Proof}: Proved in \cite{BL} (see Corollary 9.9).\hfill$\square$
\\

From now on, we assume that $G$ = $G_\Gamma$ is a right-angled Artin group of rank $r$ ($r$ $\geq$ 1). We shall show that $G$ satisfies the conditions of Theorem \ref{7.8}. Since $B_1$ is precisely the Torelli group of $G$, Theorem \ref{7.14} will then result from the following:

\begin{theorem}\label{7.9}
Let $B$ be a group. Suppose that $B$ admits a separating $\mathbbm{Z}$-linear central filtration $(B_n)_{n \geq 1}$. Then $B$ is residually torsion-free nilpotent.
\end{theorem}

\textit{Proof}: Proved in \cite{BL} (see Theorem 6.1).\hfill$\square$
\\

We need to introduce the following notations. Let $K$ be a commutative ring. We denote by $M_\Gamma$ the monoid defined by the presentation:
\begin{center}
$M_\Gamma$ = $\langle$ $V$ $\mid$ $vw$ = $wv$, $\forall$ \{$v$,$w$\} $\in$ $E$ $\rangle$,
\end{center}
by $A_\Gamma$ the associative $K$-algebra of the monoid $M_\Gamma$, and by $L_\Gamma$ the Lie $K$-algebra defined by the presentation:
\begin{center}
$L_\Gamma$ = $\langle$ $V$ $\mid$ $[v,w]$ = 0, $\forall$ \{$v$,$w$\} $\in$ $E$ $\rangle$.
\end{center}

The following theorem is due to Duchamp and Krob (see \cite{DK1}, Corollary II.16):

\begin{theorem}\label{7.10}
The $K$-module $L_\Gamma$ is free.
\end{theorem}

Thus, by the Poincar\'e-Birkhoff-Witt theorem, $L_\Gamma$ can be regarded as a Lie subalgebra of its enveloping algebra, for which Duchamp and Krob established the following (see \cite{DK1}, Corollary I.2):

\begin{theorem}\label{7.11}
The enveloping algebra of $L_\Gamma$ is isomorphic to $A_\Gamma$.
\end{theorem}

Furthermore, Duchamp and Krob proved the following (see \cite{DK2}, Theorem 2.1), which generalizes a well-known theorem of Magnus:

\begin{theorem}\label{7.12}
Suppose that $K$ = $\mathbbm{Z}$. The graded Lie $\mathbbm{Z}$-algebra associated to the lower central series of $G$ is isomorphic to $L_\Gamma$.
\end{theorem}

Set $Z$= $\cap_{v \in V}star(v)$. It follows from Servatius' Centralizer Theorem (see \cite{S}, Theorem 1) that the center $Z(G)$ of $G$ is the special subgroup of $G$ generated by $Z$.

\begin{lemma}\label{7.13}
Suppose that $Z(G)$ = \{1\}. Then $Z(L_\Gamma)$ = \{0\}.
\end{lemma}

\textit{Proof}: Let $g$ $\in$ $Z(L_\Gamma)$. Suppose that $g$ $\neq$ 0. Let $v$ $\in$ $V$. We have $[g,v]$ = 0 (in $L_\Gamma$). Now, $L_\Gamma$ can be regarded as a Lie subalgebra of $A_\Gamma$ by Theorem \ref{7.10} and Theorem \ref{7.11}. Thus, we have $gv$ = $vg$ (in $A_\Gamma$). Therefore $g$ belongs to the subalgebra of $A_\Gamma$ generated by $star(v)$ (see \cite{KR}, Theorem 2). Since $v$ is arbitrary, this leads to a contradiction with our assumption.\hfill$\square$
\\

From now on, we assume that $K$ = $\mathbbm{Z}$. We now turn to prove:

\begin{theorem}\label{7.14}
The Torelli group of a right-angled Artin group is residually torsion-free nilpotent.
\end{theorem}

\textit{Proof}: Let $H$ be the special subgroup of $G$ generated by $V \setminus Z$. Note that $Z(H)$ = \{1\}. We have: $G$ = $H\times Z(G)$. First, we show that $\mathcal{T}(G)$ = $\mathcal{T}(H)$. Let $\varphi$ : $T(H)$ $\rightarrow$ $T(G)$ be the homomorphism defined by:
\begin{center}
$\varphi(\alpha)(h,k)$ = $(\alpha(h),k)$
\end{center}
for all $\alpha$ $\in$ $T(H)$, $h$ $\in$ $H$, $k$ $\in$ $Z(G)$. Clearly, $\varphi$ is well-defined and injective. We shall show that $\varphi$ is surjective. Let $\beta$ $\in$ $T(G)$. For $g$ $\in$ $G$, we set $\beta(g)$ = $(\beta_1(g),\beta_2(g))$, where $\beta_1(g)$ $\in$ $H$ and $\beta_2(g)$ $\in$ $Z(G)$. Let $h$ $\in$ $H$. We denote by $\overline{h}$ the canonical image of $h$ in $H^{ab}$. Note that the canonical image of $h$ in $G^{ab}$ = $H^{ab} \times Z(G)$ is ($\overline{h}$,1). Since $\beta$ $\in$ $T(G)$, we have: ($\overline{h}$,1) = ($\overline{\beta_1(h)}$,$\beta_2(h)$), and then $\beta_2(h)$ = 1. Let $k$ $\in$ $Z(G)$. Since $\beta(k)$ lies in the center of $G$, we have $\beta_1(k)$ = 1. Note that the canonical image of $k$ in $G^{ab}$ is (1,$k$). As $\beta$ $\in$ $T(G)$, we have $\beta_2(k)$ = $k$. Finally, we have:
\begin{center}
$\beta(h,k)$ = $(\beta_1(h),k)$,
\end{center}
for all $h$ $\in$ $H$ and $k$ $\in$ $Z(G)$. Applying the same argument to $\beta^{-1}$, we obtain that the restriction $\alpha$ of $\beta_1$ to $H$ is an automorphism of $H$. Therefore $\beta$ = $\varphi(\alpha)$. We deduce that $\varphi$ is an isomorphism. Note that $\varphi(Inn(H))$ = $Inn(G)$. We conclude that $\mathcal{T}(G)$ = $\mathcal{T}(H)$. Thus, up to replacing $G$ by $H$, we can assume that $Z(G)$ = \{1\}. As we noted above, $\frac{G}{C^n(G)}$ is torsion-free for all $n$ $\geq$ 1. Now, for all $n$ $\geq$ 1, $C^n(G)$ $<$ $C^n_\mathbbm{Z}(G)$, and $\frac{C^n_\mathbbm{Z}(G)}{C^n(G)}$ is the torsion subgroup of $\frac{G}{C^n(G)}$ by \cite{BL}, Proposition 7.2. It follows that $C^n_\mathbbm{Z}(G)$ = $C^n(G)$ for all $n$ $\geq$ 1, and that $\mathcal{L}_\mathbbm{Z}(G)$ = $L_\Gamma$ by Theorem \ref{7.12}. Since $Z(G)$ = \{1\}, we have $Z(\mathbbm{F}_p \otimes L_\Gamma)$ = \{0\} for every prime number $p$ -- by Lemma \ref{7.13}. We deduce that $(B_n)_{n \geq 1}$ is a $\mathbbm{Z}$-linear central filtration on $\mathcal{T}(G)$ by Theorem \ref{7.8}. Now, let $\varphi$ $\in$ $Aut(G)$ such that $\varphi$ induces an inner automorphism on $\frac{G}{C^n(G)}$ for all $n$ $\geq$ 1. Let $g$ $\in$ $G$. Suppose that $\varphi(g)$ and $g$ are not conjugate in $G$. Then it follows from the proof of Theorem \ref{7.1} that there exists $n$ $\geq$ 1 such that the canonical images of $\varphi(g)$ and $g$ in $\frac{G}{C^n(G)}$ are not conjugate in $\frac{G}{C^n(G)}$ -- contradicting our assumption. Thus $\varphi$ is conjugating. Therefore $\varphi$ is inner by Proposition \ref{7.2}. We deduce that $(B_n)_{n \geq 1}$ is separating by Theorem \ref{7.8}. We conclude that $\mathcal{T}(G)$ is residually torsion-free nilpotent by Theorem \ref{7.9}.\hfill$\square$

\begin{corollary}\label{7.15}
The Torelli group of a right-angled Artin group is residually $p$-finite.
\end{corollary}

\textit{Proof}: Since $\mathcal{T}(G)$ is finitely generated by \cite{D2}, Theorem B, and residually torsion-free nilpotent by Theorem \ref{7.14}, it is residually $p$-finite by \cite{G}, Theorem 2.1.\hfill$\square$
\\

It is known that residually torsion-free nilpotent groups are bi-orderable (see, for example, \cite{CKM}, Remark 2.6). Thus, Theorem \ref{7.14} immediately yields:

\begin{corollary}\label{7.16}
The Torelli group of a right-angled Artin group is bi-orderable.
\end{corollary}

\appendix

\section{Appendix}

Let $G$ be a group, and let $H$ be a subgroup of $G$. Recall that the \textit{normal core} of $H$, denoted by $H_G$, is defined to be the largest normal subgroup of $G$ that is contained in $H$, i.e. $H_G$ = $\cap_{g\in G}gHg^{-1}$. The following lemma is probably well-known, though it does not seem to be in the literature. We include a proof for completeness.

\begin{lemma}\label{A.1}
Let $G$ be a group, and let $H$ be a subgroup of $G$. Then $H$ is open in the pro-$p$ topology on $G$ if and only if $H$ is subnormal of $p$-power index.
\end{lemma}

\textit{Proof}: If $H$ is open in the pro-$p$ topology on $G$, then it contains a normal subgroup $K$ of $p$-power index in $G$. Thus $[G:H]$ is a power of $p$. As $\frac{G}{K}$ is a finite $p$-group, every subgroup of it is subnormal. Therefore $H$ is subnormal in $G$.

Conversely, if $H$ is a subnormal subgroup of $p$-power index in $G$, then $[G:H_G]$ is a power of $p$ (see, for example, \cite{AF2}, Lemma 3.3). Thus $H$ contains an open subgroup of $G$, and hence is open itself.\hfill$\square$

\begin{flushleft}
Emmanuel Toinet, \\
Institut de Math\'ematiques de Bourgogne, UMR 5584 du CNRS, Universit\'e de Bourgogne, B.P. 47870, 21078 Dijon cedex, France \\
E-mail: Emmanuel.Toinet@u-bourgogne.fr
\end{flushleft}


\begin{thebibliography}{9}
\bibitem[AF1]{AF1}
M. Aschenbrenner, S. Friedl. 3-manifold groups are virtually residually $p$. arXiv:1004.3619.
\bibitem[AF2]{AF2}
M. Aschenbrenner, S. Friedl. Residual properties of graph manifold groups, Topology Appl. 158 (2011), 1179--1191.
\bibitem[BL]{BL}
H. Bass, A. Lubotzky. Linear-central filtrations on groups. The mathematical legacy of Wilhelm Magnus: groups, geometry and special functions (Brooklyn, NY, 1992), 45--98, Contemp. Math., 169, Amer. Math. Soc., Providence, RI, 1994.
\bibitem[CZ]{CZ}
S.C. Chagas, P.A. Zalesskii. Finite index subgroups of conjugacy separable groups. Forum Math. 21 (2009), no. 2, 347--353.
\bibitem[C]{C}
R. Charney. An introduction to right-angled Artin groups. Geom. Dedicata 125 (2007), 141--158.
\bibitem[CV1]{CV1}
R. Charney, K. Vogtmann. Finiteness properties of automorphism groups of right-angled Artin groups. Bull. Lond. Math. Soc. 41 (2009), no. 1, 94--102.
\bibitem[CV2]{CV2}
R. Charney, K. Vogtmann. Subgroups and quotients of automorphism groups of RAAGS. arXiv:0909.2444v1.
\bibitem[CKM]{CKM}
J. Cimpri$\check{\textnormal{c}}$, M. Kochetov, M. Marshall. Orderings and *-orderings on cocommutative Hopf algebras, Algebras and Representation Theory 10 (2007), no. 1, p. 25--54.
\bibitem[D1]{D1}
M.B. Day. Peak reduction and finite presentations for automorphism groups of right-angled Artin groups. Geom. Topol. 13 (2009), no. 2, 817--855.
\bibitem[D2]{D2}
M.B. Day. Symplectic structures on right-angled Artin groups: between the mapping class group and the symplectic group. Geom. Topol. 13 (2009), no. 2, 857--899.
\bibitem[DK1]{DK1}
G. Duchamp, D. Krob. The free partially commutative Lie algebra: bases and ranks. Adv. Math. 95 (1992), no. 1, 92--126.
\bibitem[DK2]{DK2}
G. Duchamp, D. Krob. The lower central series of the free partially commutative group. Semigroup Forum 45 (1992), no. 3, 385--394.
\bibitem[Dy1]{Dy1}
J.L. Dyer. Separating conjugates in free-by-finite groups. J. London Math. Soc. (2) 20 (1979), no. 2, 215--221.
\bibitem[Dy2]{Dy2}
J.L. Dyer. Separating conjugates in amalgamated free products and HNN extensions. J. Austral. Math. Soc. Ser. A 29 (1980), no. 1, 35--51.
\bibitem[G]{G}
K.W. Gruenberg. Residual properties of infinite soluble groups. Proc. London Math. Soc. (3) 7 (1957), 29--62.
\bibitem[I]{I}
E.A. Ivanova. On conjugacy $p$-separability of free products of two groups with amalgamation. Math. Notes 76 (2004), no. 3--4, 465--471.
\bibitem[KS]{KS}
A. Karrass, D. Solitar. Subgroups of HNN groups and groups with one defining relation. Canad. J. Math. 23 (1971), 627--643.
\bibitem[KR]{KR}
K.H Kim, F.W. Roush. Homology of certain algebras defined by graphs. J. Pure Appl. Algebra 17 (1980), no. 2, 179--186.
\bibitem[L]{L}
M.R. Laurence. A generating set for the automorphism group of a graph group. J. London Math. Soc. (2) 52 (1995), no. 2, 318--334.
\bibitem[Lo]{Lo}
K. Lorensen. Groups with the same cohomology as their profinite completions. J. Algebra 320 (2008), no. 4, 1704--1722.
\bibitem[LS]{LS}
R.C. Lyndon, P.E. Schupp. Combinatorial group theory. Reprint of the 1977 edition. Classics in Mathematics. Springer-Verlag, Berlin, 2001.
\bibitem[MKS]{MKS}
W. Magnus, A. Karrass, D. Solitar. Combinatorial group theory. Presentations of groups in terms of generators and relations. Reprint of the 1976 second edition. Dover Publications, Inc., Mineola, NY, 2004.
\bibitem[M]{M}
A. Minasyan. Hereditary conjugacy separability of right angled Artin groups and its applications. Groups Geom. Dyn. 6 (2012), 335--388.
\bibitem[My]{My}
A.G. Myasnikov. Approximability of outer automorphism groups of free groups of finite rank. Algebra i Logika 20 (1981), no. 3, 291--299.
\bibitem[P]{P}
L. Paris. Residual $p$ properties of mapping class groups and surface groups. Trans. Amer. Math. Soc. 361 (2009), no. 5, 2487--2507.
\bibitem[RZ]{RZ}
L. Ribes, P. Zalesskii. Profinite groups. Ergebnisse der Mathematik und ihrer Grenzgebiete. 3. Folge. A Series of Modern Surveys in Mathematics, 40. Springer-Verlag, Berlin, 2000.
\bibitem[RZ2]{RZ2}
L. Ribes, P. Zalesskii. Pro-$p$ trees and applications. New horizons in pro-$p$ groups, 75--119, Progr. Math., 184, Birkhäuser Boston, Boston, MA, 2000.
\bibitem[S]{S}
G.P. Scott. An embedding theorem for groups with a free subgroup of finite index. Bull. London Math. Soc. 6 (1974), 304--306.
\bibitem[Se]{Se}
J.P. Serre. Arbres, amalgames, ${\rm SL}_{2}$. Ast\'erisque, No. 46. Soci\'et\'e Math\'ematique de France, Paris, 1977.
\bibitem[Ser]{Ser}
H. Servatius. Automorphisms of graph groups. J. Algebra 126 (1989), no. 1, 34--60.
\end{thebibliography}
\end{document}